\documentclass{article} 
\usepackage{graphicx}
\graphicspath{{images/}}
\usepackage{amsmath,amssymb}
\usepackage{amsthm}
\usepackage{mathtools}
\usepackage{tikz}
\usepackage{subfig}
\usepackage[numbers,sort&compress]{natbib}
\usepackage[margin=1in]{geometry}
\usepackage{url}
\newtheorem{proposition}{Proposition}

\theoremstyle{plain}

\newtheorem{definition}{Definition}
\newtheorem{remark}{Remark}

\newcommand{\R}{\mathbb R}
\renewcommand{\S}{\mathbb S}
\renewcommand{\P}{\mathbb P}
\newcommand{\N}{\mathbb N}

\newcommand{\D}{\mathcal{D}}
\newcommand{\I}{\mathcal{I}}
\newcommand{\ZZ}{\mathcal{Z}}
\newcommand{\fX}{\mathfrak{X}}
\newcommand{\ind}{\mathrm{ind}}
\newcommand{\B}{\mathrm{B}}
\newcommand{\J}{\mathrm{J}}
\newcommand{\X}{\mathbf{X}}
\newcommand{\Y}{\mathbf{Y}}
\newcommand{\M}{\mathbf{M}}
\renewcommand{\L}{\mathbf{L}}

\begin{document}

\title{A bisector line field approach to interpolation of orientation fields
}


\author{Nicolas Boizot$^{\star}$  and 
        Ludovic Sacchelli$^{\dagger}$ 
}


\date{ \small $^{\star}$  Laboratoire d'Informatique et des Systèmes, Universit\'e de Toulon, Aix Marseille Univ, CNRS UMR 7020, LIS, Marseille, France\\
{\tt\small boizot@univ-tln.fr}; \\
\small $^{\dagger}$  Department of Mathematics, Lehigh University, Bethlehem, PA, USA. \\
{\tt\small lus219@lehigh.edu}; \\
~\\
This research has been partially supported by the ANR SRGI (reference ANR-15-CE40-0018). \\
This is a post-peer-review, pre-copyedit version of an article published in Journal of Mathematical Imaging and Vision (JIMIV). The final authenticated version is available online at: \url{http://dx.doi.org/10.1007/s10851-020-00990-5}
}


\maketitle

\begin{abstract}
We propose an approach to the problem of global reconstruction of an orientation field. The method is based on a geometric model called \emph{bisector line fields}, which maps a pair of vector fields to an orientation field, effectively generalizing the notion of doubling phase vector fields. Endowed with a well chosen energy minimization problem, we provide a polynomial interpolation of a target orientation field while bypassing the doubling phase step. The procedure is then illustrated with examples from fingerprint analysis.
\emph{Orientation fields \and bisector line fields \and polynomial interpolation \and fingerprint analysis \and singularities}
\end{abstract}

\section{Introduction} \label{sec:intro}

The present article deals with the question of global reconstruction of orientation fields on the basis of a discrete dataset.
The aim is to present an alternative way of modeling orientation fields that allows to 
use a natural energy.

{ As continuous mathematical objects, orientation fields adequately model texture patterns predominantly displaying orientation information. They provide a  unifying framework for various patterns observed in nature such as fingerprints \cite{Mal-Mai2009,Kas-Wit1987,Wan-Hu2007a,Gup-Tiw2019}, liquid crystals arrangements in their nematic phase \cite{Cha1992,Dec-Red2015,Deg-Pro1995,Li-Shi2019} or the pinwheel structure of the visual cortex V1 of mammals \cite{pet2017,Hub-Wie1962,Bos-Che2014,Bos-Dup2012,Cit-Sar2006,Pra-Gau2018}.}


The model we want to discuss ties in with classical techniques used in the field of fingerprint reconstruction and authentication, therefore, the problem of the estimation of fingerprint ridge topologies is used to illustrate this approach. Indeed, as it is emphasized in \cite{Wan-Hu2007a,Mal-Mai2009,Wan-Hu2011,Jos-Dey2009,Gup-Tiw2019,Kho-Kho2017}, the estimation of fingerprint ridge topologies can be a necessary step before the use of high-level classification algorithms. In the
 review article \cite{Bia-Xu2019}, the authors proposed a classification of estimation methods into three broad categories: gradient-based methods, mathematical model-based methods, and learning-based methods. The matter discussed in the present paper falls into the second one, in particular, in the sub-class of methods that don't require prior heuristic knowledge.

The classical procedure consists in first obtaining a coarse estimation of the orientation field in the form of a discrete dataset. This  step can be achieved with gradient-based methods. Then modelling choices are made and an optimisation algorithm is applied to fit a model to the dataset. Finally, the orientation field is reconstructed with the help of the fully identified mathematical model.

On the first hand, large discrete datasets can be constructed by means of a rough method such as the computation of the gradients of the fingerprint image gray intensity changes \cite{Baz-Ger2002,Kas-Wit1987}. However, such a strategy is prone to introduce significant noise into the data \cite{Wan-Hu2007a}. Note that enhancement techniques can be applied in order to improve the situation \cite{Bia-Xu2019,Bia-Din2017,Zhu-Yin2006}. 
 On the other hand, smaller but more reliable datasets can be obtained by using techniques focusing on the detection of stable and highly distinctive fingerprints features, such as minutiae \cite{Li-Kot2012,Jos-Dey2009}. In this case, the difficulty lies in the fact that the global reconstruction of the orientation field has to be performed on the basis of scarcer information.  

From the modeling point of view, following \cite{She-Mon1993,God1983,Bos-Sac2016}, an orientation field can be regarded as a mapping from $\mathbb{R}^2$ to the orientation set ---that is, the interval $[0,\pi]$ where $\pi$ is equivalent to $0$--- which complicates the parametrization task. Indeed, consider the case of vector fields, which are mappings from $\mathbb{R}^2$ to $\mathbb{R}^2$. For a sufficiently regular vector field, singularities are identified by a rather simple criteria that can be computationally handled ---{\it i.e.} points where the vector field vanishes. Meanwhile, singularities of orientation fields correspond to discontinuities which, for the more classic ones, translate into points where all  orientations accumulate. 
As a consequence, in order to efficiently handle orientation fields, one needs to propose a model that re-introduces some continuity property. One of the most popular solutions to this problem, which we call  \emph{doubling phase step}, sends each dataset angle from the orientation set to $\R^2$ by doubling its value and taking its \emph{sine} and \emph{cosine}, {\it e.g.} \cite{Kas-Wit1987,Wan-Hu2007a,Li-Kot2012,Wan-Hu2011,Yau-Li2004,Gup-Gup2016,Ram-Bischof2010,Got-Mih2009,Bia-Luo2014,Jir-Hou2011}. This resolves the orientation set's cyclicity, however, the discontinuity issue remains.

Starting from this basis, several models have been proposed in the literature, in order to reconstruct the doubling phase vector field and consequently, the orientation field \cite{She-Mon1994,Zhou-Gu2004a,Zhou-Gu2004b,Dass2004,Viz-Ger1996,Cap-Mai2007,Ros-Sha2007,Feng-Jain2011,Bia-Xu2019}. For instance, in the seminal paper \cite{Wan-Hu2007a}, the authors propose a method based on 2D Fourier Expansions to interpolate the doubling phase field, which helps address the discontinuity of the target without prior information. However, one can remark the two following points: first, this approach does not fully resolves the discontinuity issues, and second, the targeted field is not  the true orientation field.

The present paper proposes to address the problem of global reconstruction of orientation fields by means of bisector line fields. This model, discussed in \cite{Bos-Sac2016} in the framework of differential geometry,  appears as a natural extension of the doubling phase step. 

 The \emph{bisector line field} is an orientation field constructed from  two vector fields according to the following procedure. At each point in $\mathbb{R}^2$,  the two vector fields define two directions (in $[0,2\pi]$). The mean value of these two directions, taken modulo $\pi$, belongs to the orientation set and corresponds to the orientation of the line bisecting the angle between the two vector fields.  This concept displays many properties that makes it an adequate tool for the global reconstruction task in the sense that, for sufficiently regular generating vector fields, the behavior of the bisector line field is very tractable. For instance, in the practical cases considered in this paper, singularities of bisector line fields happen whenever one of the two generating vector fields vanishes. As a consequence, this model organically solves the discontinuity issue discussed earlier.

The second contribution of this article is to rely on
an energy functional that measures the error between two orientation fields directly in the orientation set. In fact, it is very similar to the \emph{Root Mean Square Deviation} for orientation fields, which is acknowledged as well suited in order to measure the distance between the original ({\it i.e.} the ground truth) and the reconstructed orientation field \cite{Tur-Mal2011,Bia-Xu2019}.
In other words, the doubling phase step is discarded, and the reconstruction task is performed by minimising a classic comparison index.
Moreover, when this energy is associated with the bisector line field model, the resulting optimization problem is particularly suited to the gradient method since the complexity of the gradient only depends on the underlying vector field structure. Therefore, although polynomial interpolation is performed in the present paper, the proposed approach might also be used with the techniques presented in \cite{Wan-Hu2007a,Ram-Bischof2010,Tas-Hel2010,Liu-Liu2014,Gup-Gup2016,Gup-Gup2015}. 

The rest of the article is organised as follows. Section~\ref{sec:bisector} is dedicated to the definition of the bisector line field model. A general presentation of the proposed approach for orientation field reconstruction is also exposed. Important mathematical properties that justifies the approach are presented in Section~\ref{sec:foundations}. It can be skipped by readers not interested in these aspects. Section~\ref{sec:energy} further details the interpolation method with the definition of the energy function and the calculation of the gradient. Finally, results are displayed and discussed in Section~\ref{sec:discussion}.

\section{Strategy overview}\label{sec:bisector}

\begin{figure}
\begin{center}
\subfloat[Fingerprint]{\includegraphics[width = 0.3 \linewidth]{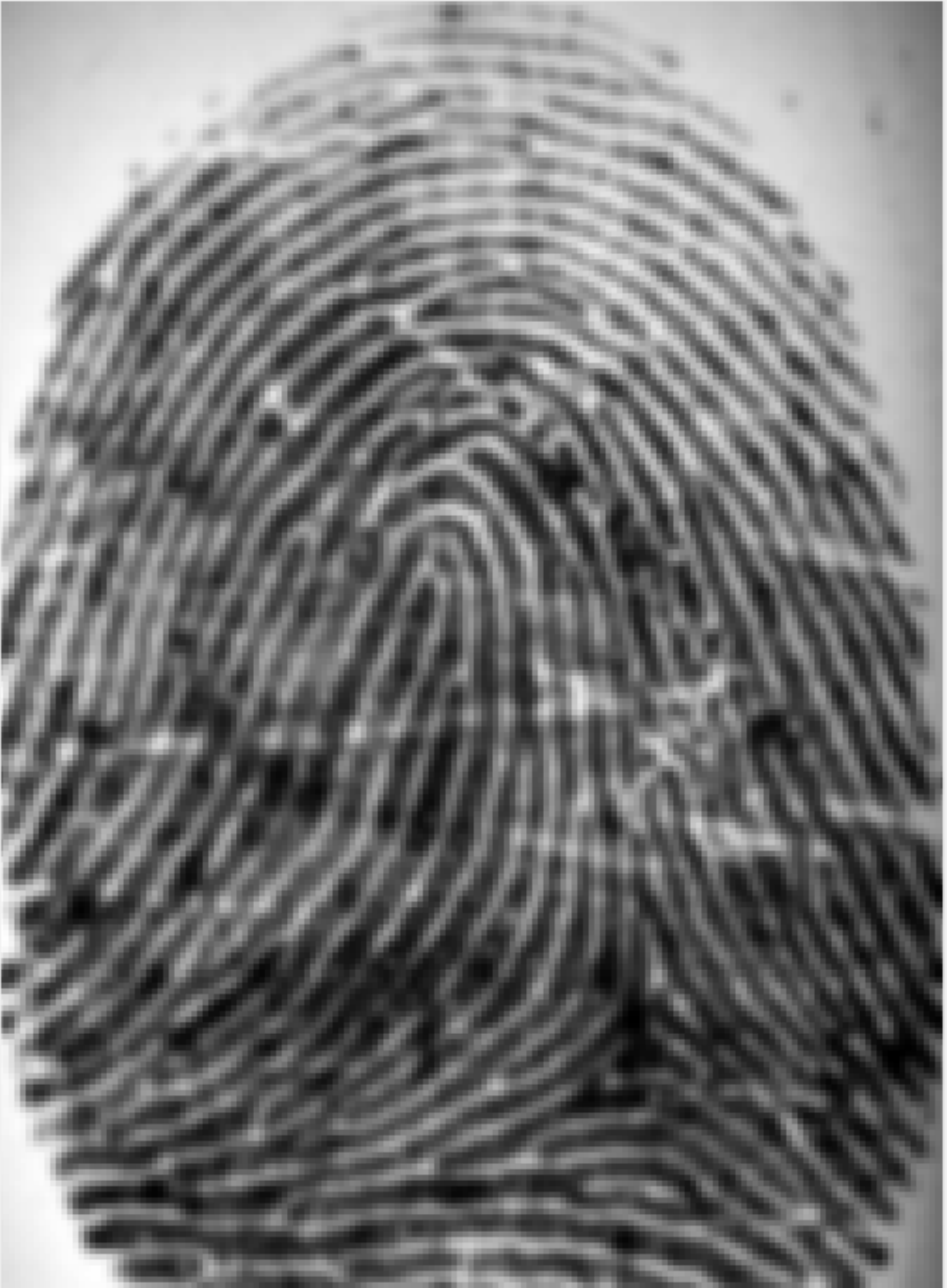} \label{fig:ReconExpl_Fingerprint}} \quad
\subfloat[Estimated OF]{\includegraphics[width = 0.3 \linewidth]{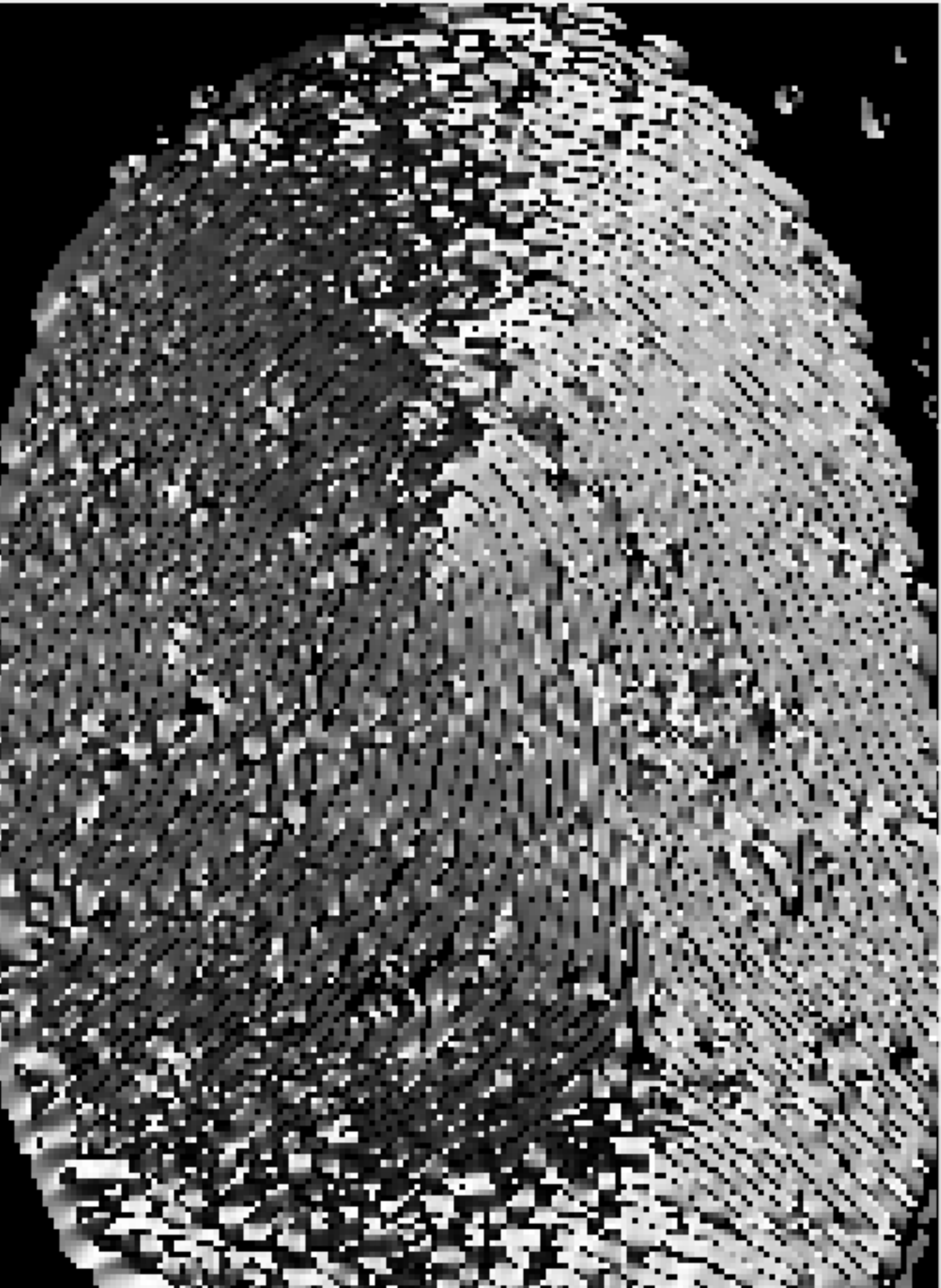}} \quad
\subfloat[Reconstructed OF]{\includegraphics[width = 0.3 \linewidth]{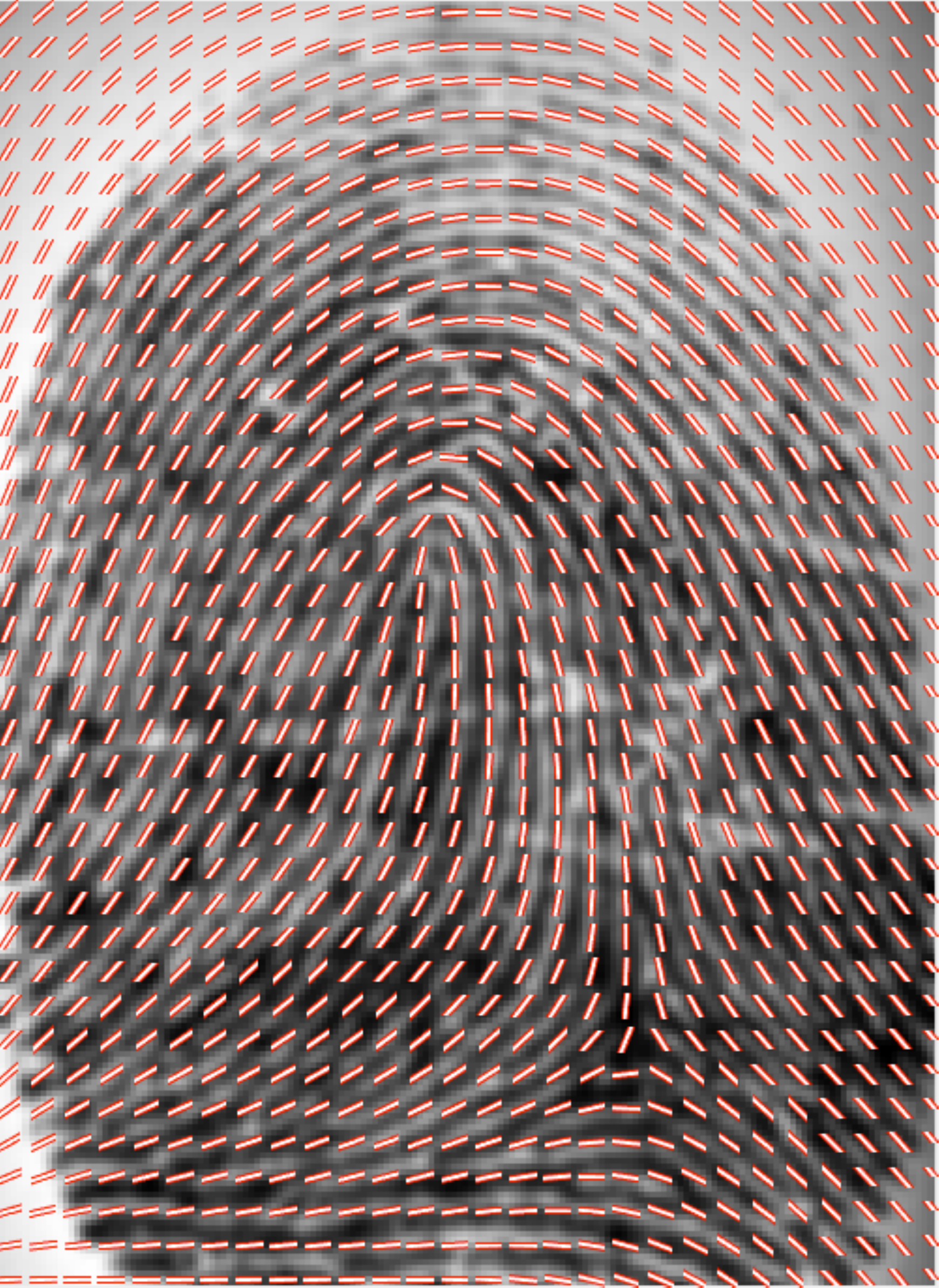}}
\caption{Reconstruction example. In (b) the grey level scale is such that black pixels correspond to $0$ and white ones to $\pi$.}
\label{fig:ReconExpl}
\end{center}
\end{figure}

The goal of orientation field interpolation is to transform a set of discretely localized orientation information into a global model, as illustrated in Figure~\ref{fig:ReconExpl}. Starting from a fingerprint image, a local estimation of the target orientation field is obtained. Enhancement techniques that improve quality of the dataset, {\it e.g.} \cite{Bia-Xu2019,Bia-Din2017,Zhu-Yin2006, Gup-Gup2016,Zhou-Gu2004a}, can  classically be used 
at this point.

Let us stress that the amount of available data can be rather small in some cases, for instance in minutiae-based extraction \cite{Ros-Sha2007} or in low quality fingerprint impression.

\begin{definition}Let  $\D$ denote a subset of $\R^2$.

\begin{itemize}
    \item  $ \S^{n}=\left\{ x \in \mathbb{R}^{n+1} \mid
\sum_{k=1}^{n+1} x_k^2 = 1   \right\}$ denotes the \emph{unit n-sphere}.
In particular $\S^1$ is interpreted as the interval $[0,2\pi]$ where $0$ is identified with $2\pi$.
    \item The \emph{orientation set}, denoted by $\P^1$ is the interval $[0,\pi]$ where $0$ is identified with $\pi$.
    \item  A \emph{vector field} is a map $\X:\D\to \R^2$. We denote its \emph{zero set} by $\ZZ_\X = \left\{ p \in \D:\X(p) = (0,0) \right \}$.
    \item  An \emph{orientation field} is a map $\L:\D\to \P^1$.
\end{itemize}
\end{definition}

\begin{figure}[htb]
\begin{center}
\includegraphics[width=0.65\columnwidth]{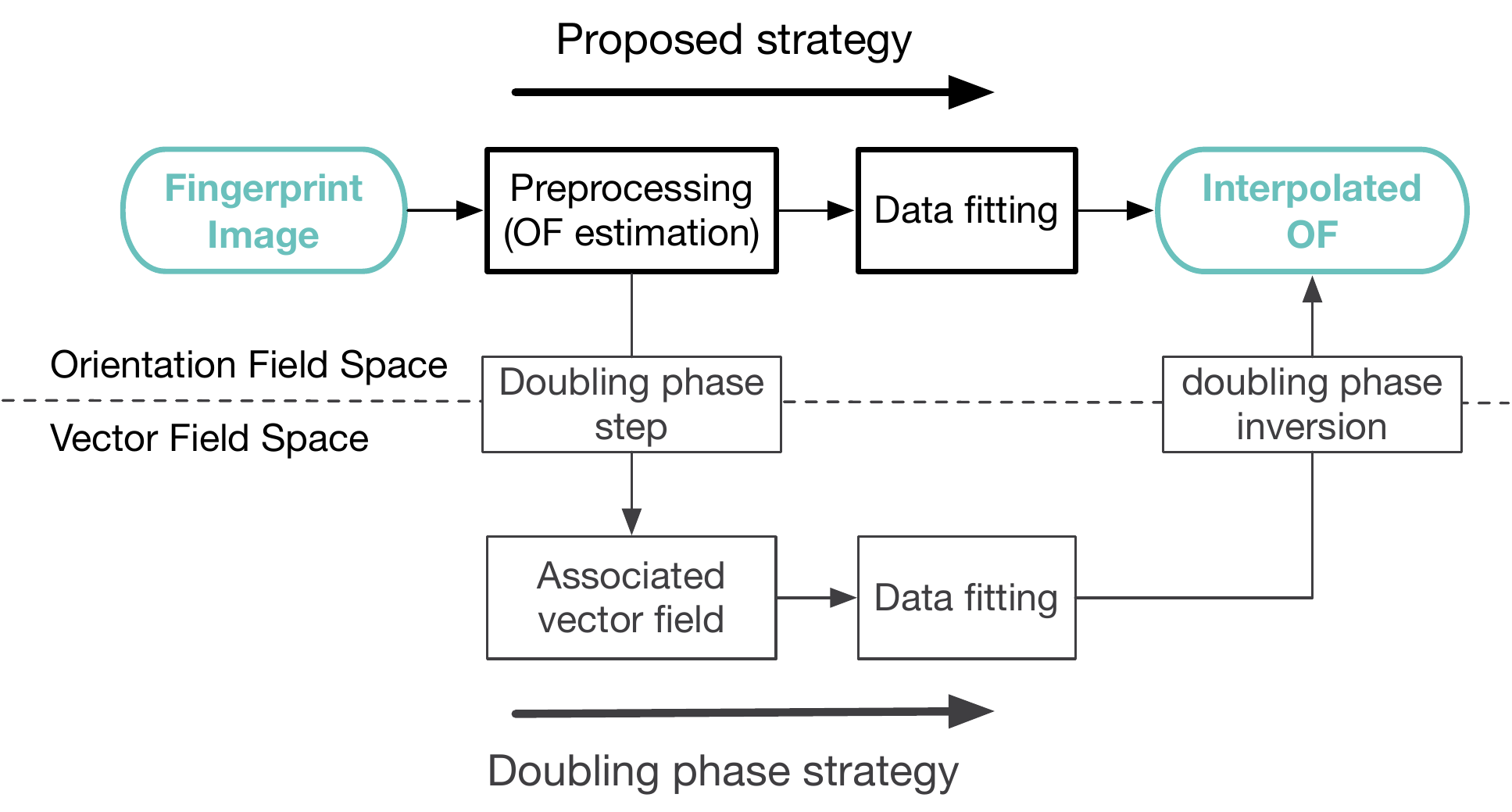} 
\caption{FlowGraph of the proposed approach}\label{fig:TheStrategy}
\end{center}
\end{figure}

A classical strategy, \cite{Kas-Wit1987,Wan-Hu2007a,Li-Kot2012,Wan-Hu2011,Yau-Li2004,Gup-Gup2016,Ram-Bischof2010,Got-Mih2009,Bia-Luo2014,Jir-Hou2011},  corresponding to the bottom branch  in  Figure~\ref{fig:TheStrategy},  is to not treat the estimation of the fingerprint orientation field as such.
Let $\left ( \theta_i \right)$ denote a collection of elements of $\P^1$, then a collection of elements $\left ( V_i \right)$ of $\R^2$ is obtained by:
$$ V_i = \begin{pmatrix} \cos{  2 \theta_i} \\ \sin{ 2 \theta_i  } \end{pmatrix} .$$
Next, data are fitted. This requires choosing both a mathematical model for the global vector field and a fitting technique. Finally, the interpolated data is taken back into the orientation set. 

As emphasized in Figure~\ref{fig:TheStrategy}, we propose a more straightforward strategy based on a model called the {\it bisector line field}.

Similarly to vector fields, orientation fields admit an intrinsic definition from the differential topology point of view, see \cite{Bos-Sac2016} for a thorough discussion. However, such a definition lacks efficiency if one hopes to work with orientation fields and apply techniques from differential calculus. The introduction of bisector line fields is meant to bridge this gap.
%
%
\begin{definition}
Let  $\D$ denote a subset of $\R^2$.  
Let $\X$ and $\Y$ be two vector fields on $\D$. 
The \emph{bisector line field of $(\X,\Y)$} is the orientation field
$$
\B(\X,\Y):\D\setminus(\ZZ_\X\cup \ZZ_\Y)\to \P^1
$$
 such that at any point $p\in\D$, $\B(\X,\Y)(p)$ is the orientation of the line bisecting the oriented pair of vectors $(\X(p),\Y(p))$.
\end{definition}

\begin{figure}[htb]
\begin{center}
\begin{tikzpicture}[line join=round,>= latex,x=2.0cm,y=2.0cm]

\clip(1.4,1.8) rectangle (5.,3.8);

\draw [->,line width=1.pt,color=black] (2.8326419725368748,2.2164129218750697) -- (4.685595354668596,3.420621357556255);

\draw [->,line width=1.pt,color=black] (2.8326419725368748,2.2164129218750697) -- (1.498242498948219,3.359383032014635);

\draw [line width=1.pt,domain=1.4:5.,color=black] plot(\x,{(-2.6803154837669454--0.9978232810591935*\x)/0.0659446720839969});

%

\draw [shift={(2.8326419725368748,2.2164129218750697)},->,line width=1.pt,color=black] (0.5553920758325129:0.6) arc (0.5553920758325129:86.21890474934047:0.6);

%

\draw [->,line width=1.pt,color=black] (2.8326419725368748,2.2164129218750697) -- (3.832594991705581,2.2261061904357087);

\draw [shift={(2.8326419725368748,2.2164129218750697)},->,line width=1.pt,color=black] (0.5553920758325129:0.8) arc (0.5553920758325129:33.019275032959115:0.8);

\draw [shift={(2.8326419725368748,2.2164129218750697)},->,line width=1.pt,color=black] (0.5553920758325129:0.4) arc
(0.5553920758325129:139.4185344657218:0.4) ;

\draw[color=black] (4.7,3.2) node {$\X$};
\draw[color=black] (1.65,3.45) node {$\Y$};
\draw[color=black] (3.8,2.05) node {$\mathbf{Z}$};
\draw[color=black] (3.3,3.6) node {$\B(\X,\Y)$};
\draw[color=black] (3.25,2.9) node {$\frac{\theta_\X+\theta_\Y}{2}$};
\draw[color=black] (3.8,2.5) node {$\theta_\X$};
\draw[color=black] (2.65,2.75) node {$\theta_\Y$};

\end{tikzpicture}
\caption{ Representation of the bisection operation on an oriented pair of vectors, where $\mathbf{Z}=(1,0)'$ marks the reference vector to measure angles.}
\label{F:bissection}
\end{center}
\end{figure}
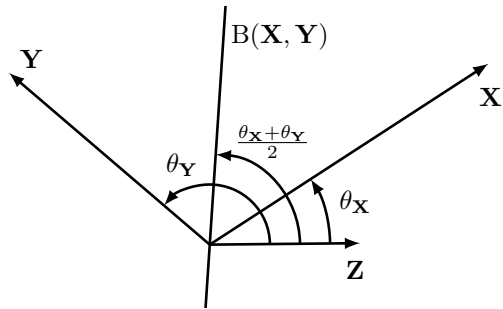

In practice, this definition implies the following.
Let $p\in \D\setminus(\ZZ_\X\cup \ZZ_\Y)$, there exist unique $\theta_\X(p)$ and $\theta_\Y(p)$ in $\S^1$ such that 
$$
\X(p)=\|\X(p)\|
\begin{pmatrix}
\cos\theta_\X(p)
\\
\sin\theta_\X(p)
\end{pmatrix}$$
and
$$
\Y(p)=\|\Y(p)\|
\begin{pmatrix}
\cos\theta_\Y(p)
\\
\sin\theta_\Y(p)
\end{pmatrix}.
$$
Then, as illustrated in  Figure~\ref{F:bissection},
$$
\B(\X,\Y)(p)=\frac{1}{2} \left(\theta_\X(p)+\theta_\Y(p)\right)\pmod{\pi}.
$$


Our strategy consists in building $\X$ and $\Y$ such that $B\left(\X,\Y  \right)$ interpolates the data. The data fitting step then relies on the above formula. In this regard, the modelling of the vector fields $\X$ and $\Y$ is completely open. Numerous strategies have been identified for this purpose, be it  rational complex functions \cite{Zhou-Gu2004b}, Taylor expansions \cite{For-Str1993,Yau-Li2004}, orthogonal/Legendre/Chebychev polynomials \cite{Tas-Hel2010,Liu-Liu2014,Ram-Bischof2010,Bia-Luo2014}, Fourier expansions \cite{Wan-Hu2007a,Tas-Hel2009,Tao-Yan2010}, and others \cite{Bia-Xu2019}. For the sake of simplicity, we chose to implement a Taylor expansion model.

\begin{remark}
As described in the introduction, the concept of orientation field arises in different fields of applications. However, as it is a rather uncommon notion,  terminology can vary between fields of study. In the context of differential geometry, the term line field is sometimes preferred, following Hopf's terminology ``fields of line elements'' \cite{Hopf,Bos-Sac2016}. We chose the terminology that is more popular in the fingerprint analysis  community.
\end{remark}

\section{Bisector line fields}\label{sec:foundations}

{
In this section, we expose properties of  bisector line fields  that highlight the suitability of this structure for the  interpolation of orientation fields.
}
\subsection{Singularities of line fields}

In our application case, singularities of orientation fields should be understood as unresolvable discontinuities of the map $\L:\R^2\to \P^1$. To this extent, we refer in the following to continuous orientation fields $\L:\D\to \P^1$, where $\D$ is the continuity domain of $\L$, and singularities of $\L$ are located in $\bar{\D}\setminus \D$, the border of $\D$, where $\bar{\D}$ denotes the topological closure of $\D$.

\begin{definition}
A continuous orientation field $\L:\D\to \P^1$ is said to be \emph{singular} at a point $p\in \bar{\D}\setminus \D$ if it cannot be uniquely continuously extended at $p$.
\end{definition}

In the following, we focus on isolated singularities, which are sufficient to our framework.
 In particular, we can introduce the topological index, which is a useful tool for the study of such singularities of line fields.

Let $\L:\D\to \P^1$ be a continuous orientation field and let $p$ be an isolated singularity of $\L$ such that $p$ belongs to the interior of $\bar{\D}$ (for instance there exists an open subset $V$ of $\R^2$   such that $V=\{p\}\cup(V\cap\D)$). Then the index of $V$ at $p$, denoted $\ind_p(\L)$, is a half integer quantifying the winding of $\L$ around $p$  (for a precise construction of this object see, for instance, \cite[Section 3.2]{Bos-Sac2016}).

In the context of the study of orientation fields singularities, the following properties of bisector line fields arise.

\begin{proposition}[\cite{Bos-Sac2016}]\label{P:continuity}
 If $\X$ and $\Y$ are continuous vector fields over $\R^2$ then $\B(\X,\Y)$ is continuous on $\D=\R^2\setminus(\ZZ_X\cup \ZZ_\Y)$.
 
On the other hand, let $\L:\D\to \P^1$ be a $C^k$ orientation field, for some $k\in \N\cup\{\infty\}$. Then there exist two $C^k$ vector fields $\X,\Y$ over $\R^2$ such that $\L=\B(\X,\Y)$.
\end{proposition}

\begin{proposition}
Let $(\X, \Y )$ be a pair of continuous vector fields over $\R^2$. Given an isolated point $p$ of $\ZZ_\X \cup \ZZ_\Y$, we have 
$$
\ind_p(\B(\X, \Y )) =  \frac{1}{2}(\ind_p(\X) + \ind_p(\Y )).
$$
\end{proposition}

\begin{remark}
An immediate consequence of this model is that singularities, initially characterized by an analytical property, now coincide with the zero sets of regular functions.
In particular, isolated singularities of the bisector line field corresponds to  isolated singularities of either one of the two generating vector fields.
\end{remark}

\subsection{Generic properties of bisector line fields}

Following Thom \cite{thom1972stabilite}, it is understood, philosophi\-cally speaking, that ``typical'' behaviors of mathematical objects should be the only ones visible in nature. From the point of view of transversality theory, these typical properties are known as \emph{generic}, in the sense that a property is generic on a topological set if it is satisfied on a residual subset.
A residual subset is understood to be a large dense set
in the following sense:
it is a countable intersection of open and dense subsets.

The classical application case of this theory is geared towards regular maps, endowed with the $C^k$ Whitney topology, for some $k\in \N$ (see, {\it e.g.} \cite{Hir1976,Abr-Rob1967}).

When generic features of vector fields are considered, more can be said \cite{Bos-Sac2016}. Example of such features are isolation of singularities of vector fields or the fact that indices of such singularities must be  $\pm 1$.

Let $\X$ and $\Y$ be two smooth vector fields. Generically with respect to the $C^2$ Whitney topology, $\ZZ_\X$ and $\ZZ_\Y$ are discrete collections of points that do not accumulate.
Thus, singularities of $\X$ and $\Y$ are isolated, and moreover have $\pm 1$ index.
Furthermore, when considering $(\X,\Y)$ as a pair, we get  $\ZZ_\X\cap \ZZ_\Y=\emptyset$.

As a consequence of these facts and Proposition~\ref{P:continuity}, we get the following property for generic smooth bisector line fields, as illustrated in Figure~\ref{F:Bissection_example}.

\begin{proposition}
Let $(\X,\Y)$ be a generic pair of vector fields. 
Singularities of $\B(\X,\Y)$ are isolated and have index $\pm1/2$.
\end{proposition}

\begin{remark}
The set of orientation fields is not endowed with a topology allowing study of genericity or stability. It is then necessary to provide an extrinsic structure and this was the motivation for the study of bisector line fields in \cite{Bos-Sac2016}. Another popular model relies on differential 2-forms on surfaces \cite{Sot-Gar2008}.
\end{remark}

{
\begin{figure}[hbtp]
\begin{center}
\subfloat[Vector field $\X$]{
\includegraphics[width=.3\linewidth,trim=1mm 1mm 1mm 1mm, clip=true]{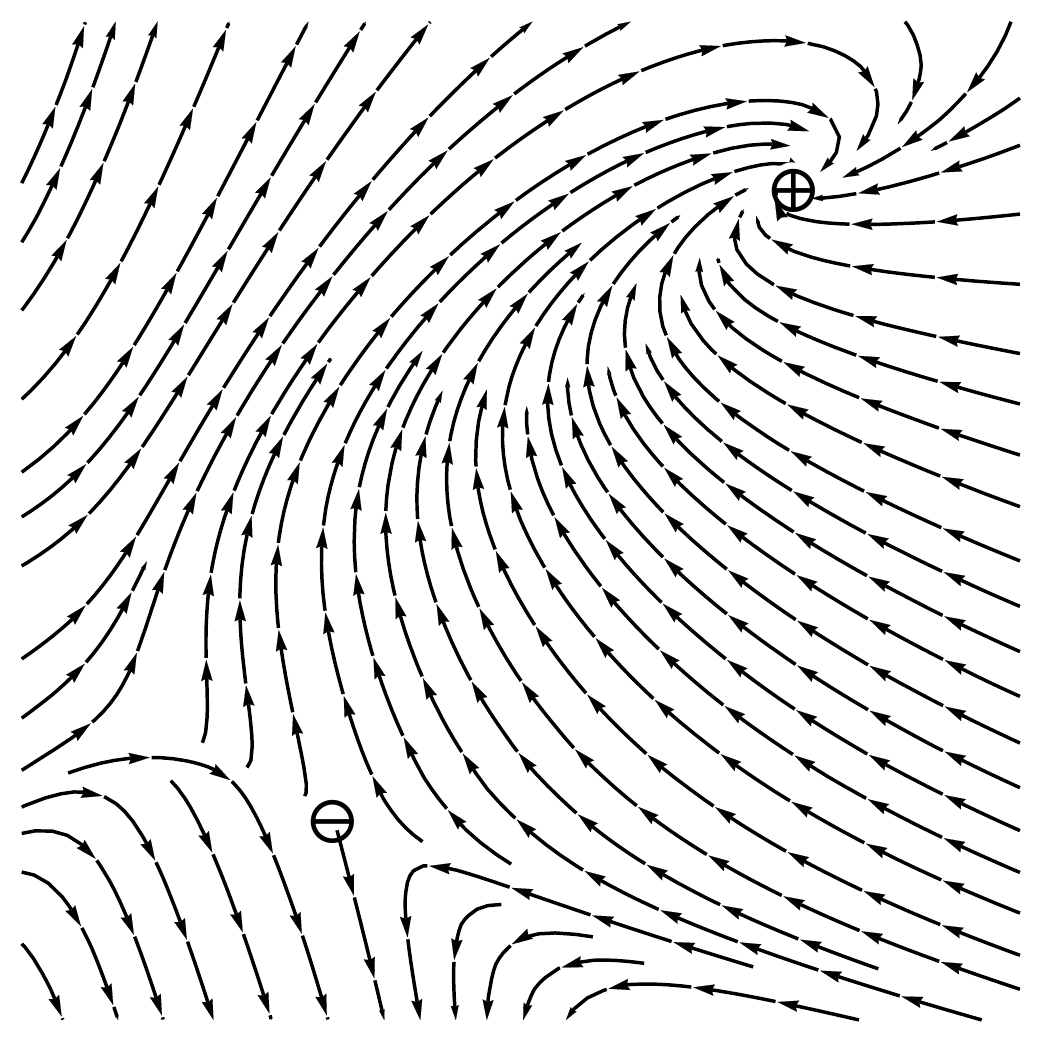}}
\quad
\subfloat[Vector field $\Y$]{
\includegraphics[width=.3\linewidth,trim=1mm 1mm 1mm 1mm, clip=true]{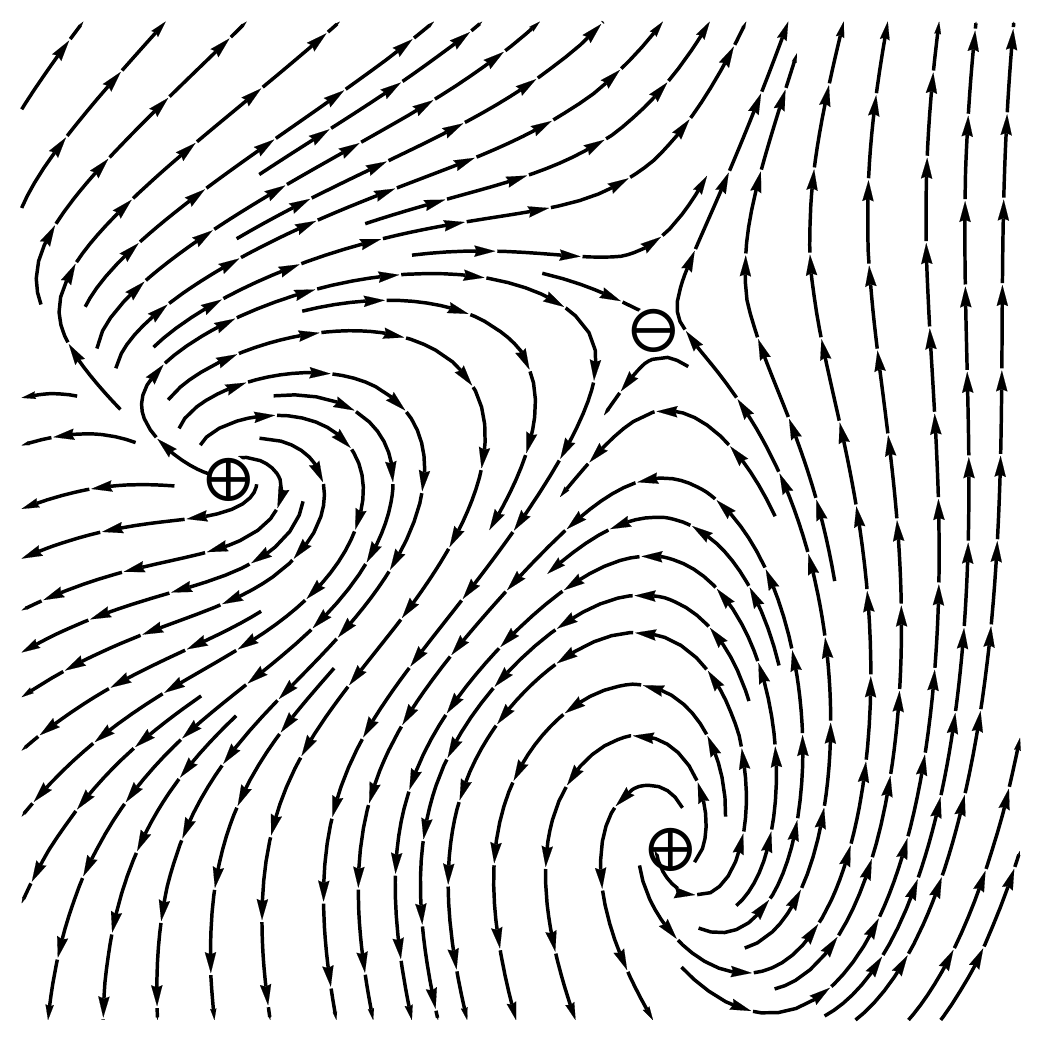}}
\quad
\subfloat[Orientation field $\B(\X,\Y)$]{
\includegraphics[width=.3\linewidth,trim=1mm 1mm 1mm 1mm, clip=true]{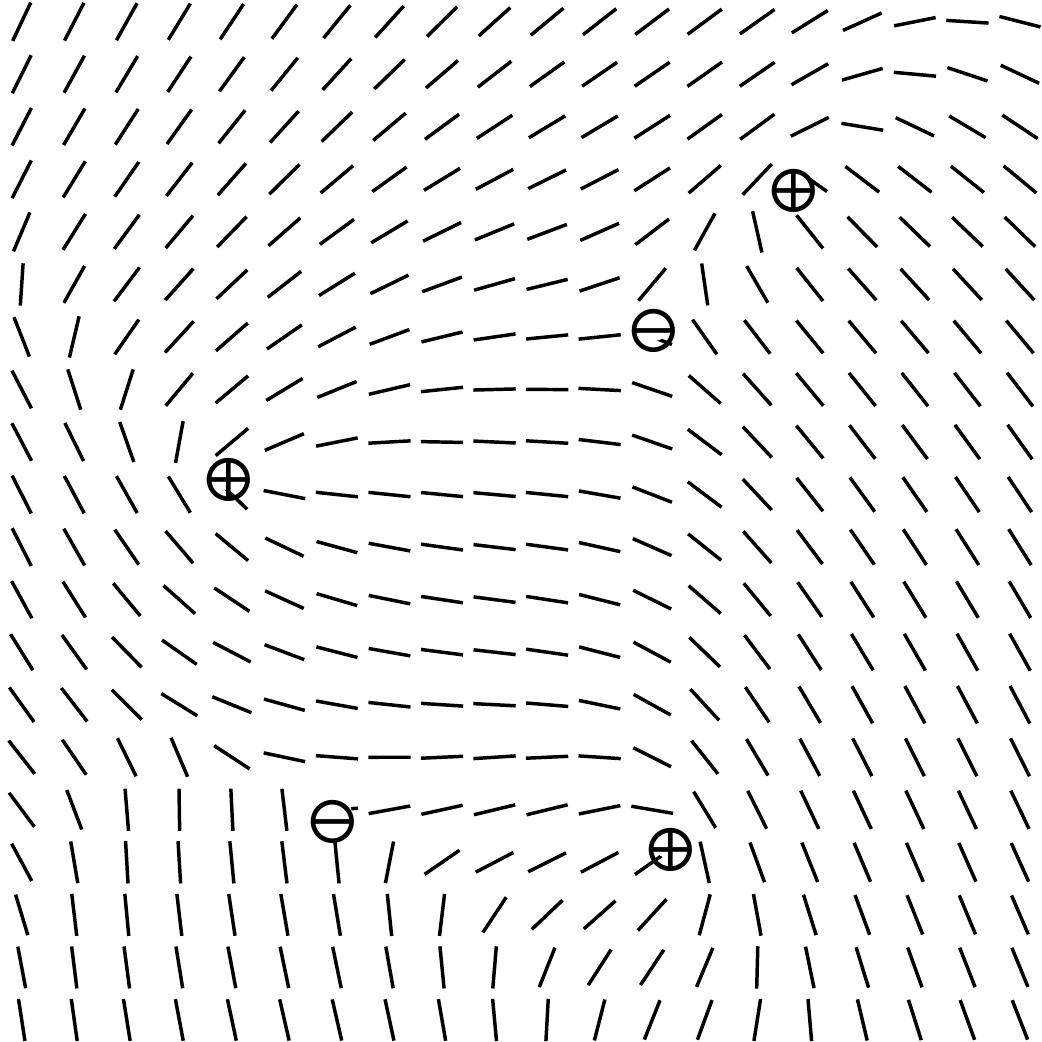}}
\end{center}
\caption{
Example of the bisection of two vector fields.
Singularities of positive index have been marked with a { $\oplus$} symbol, while singularities of negative index have been marked with a $\ominus$ symbol.
}
\label{F:Bissection_example}
\end{figure}

}

\section{Interpolation as an energy minimization problem}\label{sec:energy}

This section presents the elements needed the implementation of the orientation field reconstruction task. That is, the energy function, the polynomial model used for vector fields, and the gradient descent scheme. 

\subsection{Discrete energy}\label{subsec:DiscEner}

The \emph{Root Mean Square Deviation} for orientation fields is acknowledged as well suited to measure the distance between two orientation fields \cite{Tur-Mal2011,Bia-Xu2019}. As seen in Section~\ref{sec:bisector}, the bisector line field model allows to work in the orientation set. Therefore, the root mean square deviation can now be directly used to perform the reconstruction task. This section details this quantity in the energy formalism, which is suitable for gradient descent schemes. Let us underline that, for a given number of data, minimizing the energy is equivalent to minimizing the root mean square deviation.

Let $d:[0,\pi]\times [0,\pi]\to [-\pi/2,\pi/2] $ be piecewise defined by
\begin{equation}
 d(\theta_1,\theta_2)=\left\{ 
\begin{aligned}
 &\theta_1-\theta_2 & \text{ if }& \left|\theta_1-\theta_2\right| \leq \frac{\pi }{2}, \\
&\theta_1-\theta_2-\pi  & \text{ if }&\theta_1-\theta_2>\frac{\pi }{2} ,
\\
&\theta_1-\theta_2+\pi  & \text{ if }&\theta_1-\theta_2<-\frac{\pi }{2} .
\end{aligned}
 \right.
 \label{equ:distance}
\end{equation}
 Let $\mathcal{I} = \left \{1,\hdots,m \right \} $ where $m$ is the number of data, and let $(x_i,y_i)_{i\in \mathcal{I}}$ be a  collection of points in $\R^2$.
For $\M,\L:\D\subset \R^2\to \P^1$, such that $(x_i,y_i)_{i\in \I}\subset \D$,
we set  the least-squares energy functional to be
$$
\J(\M,\L)=\sum_{i\in \I} d\left(\M(x_i,y_i),\L(x_i,y_i)\right)^2. 
$$

The relevance of $d$ is further justified by the following observations.

\begin{figure}[hbtp]
\centering
\hfill
\subfloat[Flat representation]{
\includegraphics[width=.4\linewidth,trim=1mm 1mm 1mm 1mm, clip=true]{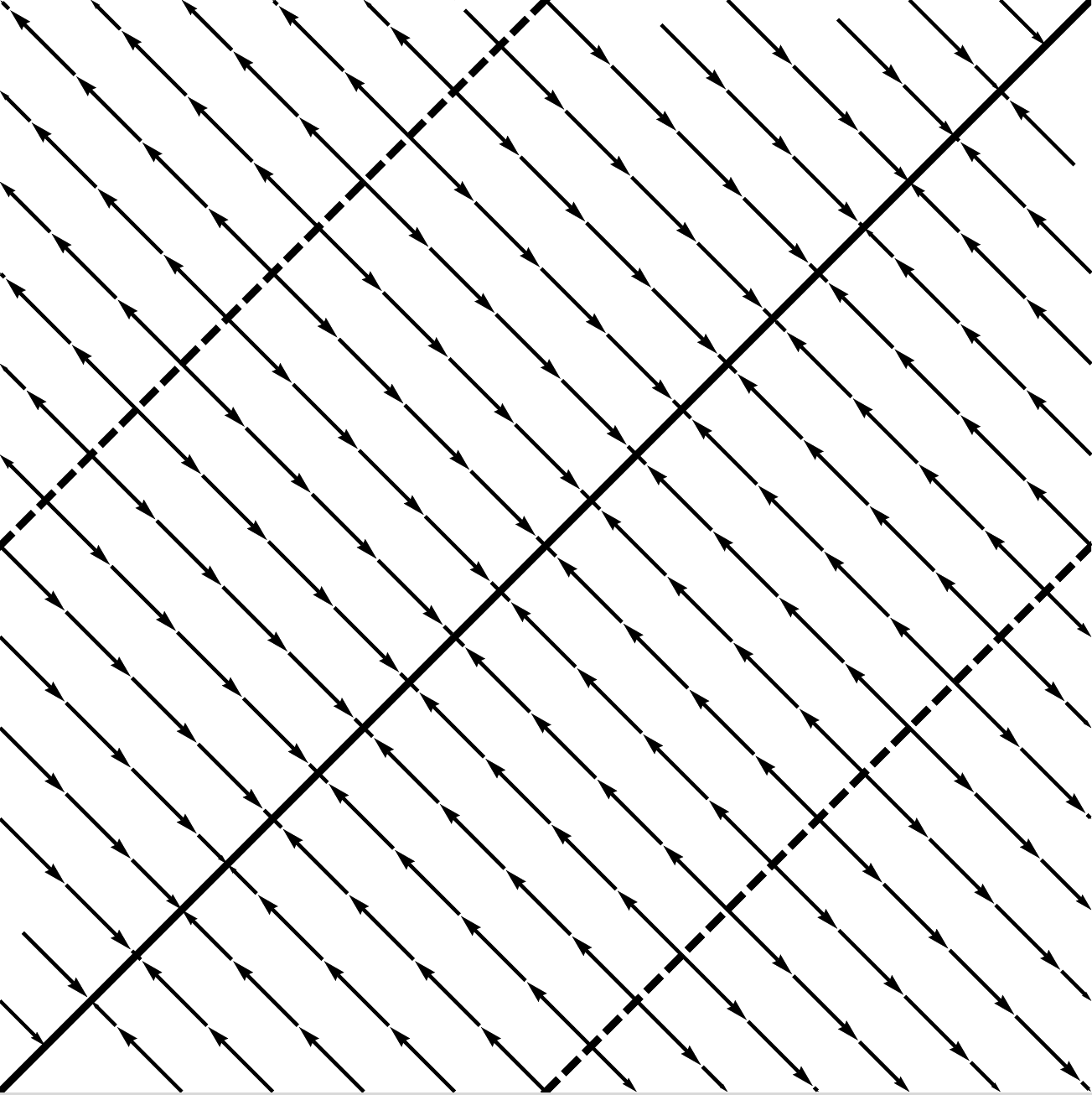}}
\hfill
\subfloat[Embedded 3d representation]{
\includegraphics[width=.4\linewidth,trim=1mm 1mm 1mm 1mm, clip=true]{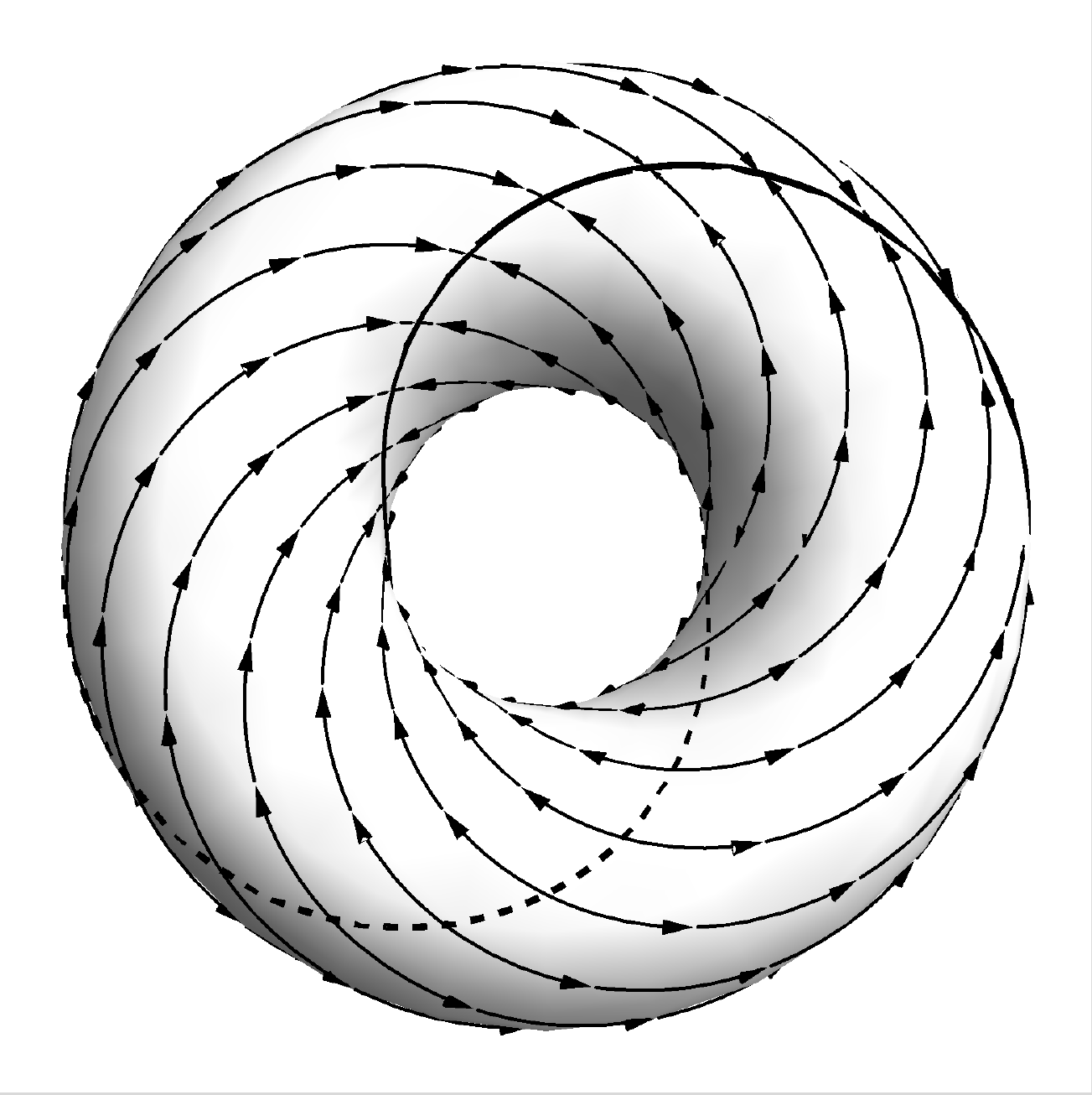}}
\hfill{}
\caption{
Representation of the gradient flow of $d^2$ on the flat torus $\P^1\times \P^1$. The set $\{\theta_1-\theta_2=\pm \pi/2\}$ is denoted with a dashed line, the set  $\{\theta_1=\theta_2\}$ is denoted with a continuous line.}
\label{F:torus}
\end{figure}

\begin{remark}\label{R:torus}

\begin{enumerate}
\item
Consider the energy given by $d^2$ on the torus $\P^1\times \P^1$. The gradient of $d^2$ is 
$$
\nabla d^2(\theta_1,\theta_2)= 2 d(\theta_1,\theta_2)(\partial_{\theta_1}-\partial_{\theta_2})
$$
As a consequence, one can check that the gradient flow of $d^2$ is actually parallel to the geodesic flow on the flat torus starting from the affine sets $\{\theta_1-\theta_2=\pm \pi/2\}$ to the diagonal $\{\theta_1=\theta_2\}$ as it is illustrated in 
 Figure~\ref{F:torus}.
 \item 
 The two connected components of $\P^1\times \P^1\setminus\{\theta_1 =\theta_2 \pmod{\pi/2}\}$ are geodesically convex once $\P^1\times \P^1$ has been endowed with the flat torus Riemannian metric. Furthermore, 
 the maps $d_{|C}$ and $d^2_{|C}$ are then geodesically convex  and strictly geodesically convex respectively.
\item Notice that $\J^{1/2}$ 
is a pseudo-metric (it is symmetric and satisfies the triangular inequality) on the space of $\P^1$-valued maps over a domain containing $(x_i,y_i)_{i\in \mathcal{I}}$.

Indeed $|d|$ is a distance over $\P^1$ and we {\it de facto} have the classical product metric over $\left(\P^1\right)^{m}$ given by
$$
d_m((\theta_1,\dots,\theta_m),(\theta_1',\dots,\theta_m'))=
\left(
\sum_{i=1}^m
d(\theta_i,\theta_i')^2
\right)^{1/2}
$$
Then with $\I = \left \{1,...,m\right\}$,
$$
(\theta_1,\dots,\theta_m)=(\M(x_1,y_1),\dots,\M(x_{m},y_{m}))
$$ and 
$$
(\theta_1',\dots,\theta_m')=(\L(x_1,y_1),\dots,\L(x_{m},y_{m})),
$$
one has
$$
\J(\M,\L)^{1/2}
=
d_m((\theta_1,\dots,\theta_m),(\theta_1',\dots,\theta_m')).
$$
\end{enumerate}
\end{remark}

\subsection{Polynomial bisector line fields} \label{subsec:PolyModel}

For a given target $\M$, we minimize $\J(\M,\L)$ over the set of polynomial bisector line fields of a fixed maximal degree.

Notice that for a given vector field $\X$, $\lambda \X$ yields the same direction $\theta_\X$ for all $\lambda>0$. Hence we can compactify the set of polynomial vector fields under consideration.
To this extent, let us define for all $n\in \N$ the set $\fX^n$ of polynomial vector fields of degree $n$ such that $\X\in \fX^n$ if there exist $(\alpha_{k,j})_{0\leq j\leq k\leq n}$, $(\beta_{k,j})_{0\leq j\leq k\leq n}$
such that 
$$
\X(x,y)=\sum_{k=0}^{n}\sum_{j=0}^k
\begin{pmatrix}
\alpha_{k,j}
\\
\beta_{k,j}
\end{pmatrix}
x^{k-j}y^j,
$$
and 
$$
\sum_{k=0}^{n}\sum_{j=0}^k
\alpha_{k,j}^2
+
\beta_{k,j}^2
=1.
$$
In particular $\left(\alpha_{k,j},\beta_{k,j}\right)_{0\leq j\leq k}\in \S^{(n+1)(n+2)}$ and the space $\fX^n$ is compact for any given $n\in \N$.

The polynomial interpolation of the target $\M$ is performed by finding $\X$ and $\Y$ solving, for $m_\X, m_\Y\in \N$,
\begin{equation}\label{E:min_problem}
    \min
    \left\{
    \J(\M,\B(\X,\Y))
    \mid
    \X\in \fX^{m_\X}, \Y\in \fX^{m_\Y}
    \right\}.
\end{equation}

Let us denote
\begin{enumerate}
 \item $n_\X =  (m_\X+1)(m_\X+2)$ and
$ n_\Y =  (m_\Y+1)(m_\Y+2)$;
    \item 
     $ \omega \in \S^{n_\X}\times \S^{n_\Y}$ such that
      $$\omega=\left(\left(\alpha_{k,j},\beta_{k,j}\right)_{0\leq j\leq k\leq m_\X},\left(\gamma_{k,j},\delta_{k,j}\right)_{0\leq j\leq k\leq m_\Y}\right);$$
    \item $\X_\omega(x,y)= \displaystyle\sum_{k=0}^{m_\X}\displaystyle\sum_{j=0}^k
\begin{pmatrix}
\alpha_{k,j}
\\
\beta_{k,j}
\end{pmatrix}
x^{k-j}y^j$

 and 
 
 $
\Y_\omega(x,y)= \displaystyle \sum_{k=0}^{m_\Y} \displaystyle\sum_{j=0}^k
\begin{pmatrix}
\gamma_{k,j}
\\
\delta_{k,j}
\end{pmatrix}
x^{k-j}y^j;
$
\item $\L_\omega=\B(\X_\omega,\Y_\omega)$.
\end{enumerate}

Problem~\eqref{E:min_problem} is then equivalent to solving
\begin{equation}
\label{E:min_problem_compact}
\min
\left\{
\J(\M,\L_\omega)\mid \omega \in  \S^{n_\X}\times \S^{n_\Y}
\right\}.
\end{equation}

The existence of minimizers is guaranteed by the compacity of $\S^{n_\X}\times \S^{n_\Y}$.

\subsection{Projected gradient descent approach}  \label{subsec:GradientDesc}
The function $d$, given in Equation~\eqref{equ:distance}, and Problem~\eqref{E:min_problem_compact} have been designed to allow the use of gradient based optimisation algorithms. Indeed, since
$$
\L_{\omega}(x,y)=\frac{1}{2}\left(\theta_{\X_\omega}(x,y)+\theta_{\Y_\omega}(x,y)\right),
$$
we have on $\R^2\setminus \ZZ_{\X_\omega}\cup \ZZ_{\Y_\omega}$
\begin{align}
\nonumber
\nabla_\omega \L_\omega
= &
\frac{1}{2}
\left(\frac{\X_{\omega,1} \nabla_\omega \X_{\omega,2}-\X_{\omega,2} \nabla_\omega \X_{\omega,1}}{\X_{\omega,1}^2+\X_{\omega,2}^2} \right.
 \\ &
 \hspace{30pt}
  \left. +
\frac{\Y_{\omega,1} \nabla_\omega \Y_{\omega,2}-\Y_{\omega,2} \nabla_\omega \Y_{\omega,1}}{\Y_{\omega,1}^2+\Y_{\omega,2}^2}
\right),
\label{E:derivative_B}
\end{align}
where $(\cdot)_{\omega,i}$ denotes the $i^{th}$ component of $(\cdot)_{\omega}$.

As a consequence, we can immediately compute by the  chain rule
\begin{multline*}
    \nabla_\omega \J(\M,\L_\omega)=\\-2\sum_{i\in \I} d\left(\M(x_i,y_i),\L_\omega(x_i,y_i)\right)
\nabla_\omega \L_\omega(x_i,y_i).
\end{multline*}
For $\rho>0$, it  should be noted  that in general:
$$ \omega+\rho \nabla_\omega \J(\M,\L_\omega)\notin \S^{n_\X}\times \S^{n_\Y}.$$ However, we solve this issue with the  projection
$$
\begin{array}{cccc}
    p :& \R^{n_\X}\times \R^{n_\Y}&\longrightarrow&\S^{n_\X}\times \S^{n_\Y} \\
    &(\omega_1,\omega_2) &\longmapsto&
    \left(
    \dfrac{\omega_1}{\|\omega_1\|}
    ,
    \dfrac{\omega_2}{\|\omega_2\|}
    \right).
\end{array}
$$

Thus in the gradient descent scheme, we iterate the recursive transform
$$
\omega_{k+1}\leftarrow p 
\left(
\omega_k+\rho \nabla_\omega \J\left(\M,\L_{\omega_k}\right)
\right).
$$

\section{Examples and discussion}\label{sec:discussion}

\subsection{Experimental setting}

In order to illustrate the paper's approach,
simulations have been performed through an implementation of a constant step-size gradient descent method in Matlab.
In line with the main application of this theory ---fingerprint analysis---  we tested the method on elements of the  FVConGoing Initiative data set  \cite{Dor-Cap2009}.

More precisely, a target has been obtained with the classical elementary method of lifting the orientations of the finger ridges from the gradient of a grey-scale image of a fingerprint \cite{Kas-Wit1987,Mal-Mai2009}. In order to discard parts of the image that don't correspond to the fingerprint, we neglected orientations 
where the gradient's norm was too low. As a consequence, 
a few data corresponding to the actual fingerprint might be missing without significant effect on the result.

In the following, we present two reconstruction examples obtained with datasets of orientations of large size but low quality ---{\it i.e.} $245 \times145-$pixels images leading to approximately 30000 data-points. One has been performed on a classical loop fingerprint, the other on a classical whorl fingerprint. A third reconstruction example has been obtained with a scarcer dataset of high fidelity, similar to the minutiae based methods found in \cite{Feng-Jain2011,Ros-Sha2007}. This reconstruction, based on a smaller dataset made of $40$ elements,  also serves to illustrate situations when latent fingerprints are treated, or the image corrupted. The scope of the present paper is not to propose a method to extract datasets in these degraded situations. However, as it is explained in Section~\ref{sec:bisector} and schematized in Figure~\ref{fig:TheStrategy}, any method that efficiently extracts data can be used prior to the bisector line field based reconstruction. 
 
Additional reconstruction experiments with datasets made of $10$, $20$, $40$ and $80$ orientations values have been performed. In each of these four cases, ten datasets where randomly generated. The results are compared to a reference set by computing the root mean square deviation \cite{Tur-Mal2011,Gup-Gup2016,Bia-Xu2019}, defined as follows.

Let $\mathcal{J}\subset \N$ and  $(\bar x_j,\bar y_j)_{j \in \mathcal{J}}$ be the collection of points ranging all the relevant pixel positions of the fingerprint image. Let  $\overline{\M}=\left(\M(\bar x_j,\bar y_j)\right)_{j \in \mathcal{J}}$ and $\overline{\L}=\left(\L(\bar x_j,\bar y_j)\right)_{j \in \mathcal{J}}$ be two matrices of orientations. Then, the root mean square deviation (RMSD) is:
 $$ RMSD (\overline{\M}, \overline{\L}) = \displaystyle \sqrt{  \frac{\sum_{j\in \mathcal{J}} d \left({\M}(\bar x_j,\bar y_j),{\L}(\bar x_j,\bar y_j)\right)^2}{\vert \mathcal{J} \vert   }  }  $$
The considered reference set is the orientation field that was reconstructed in the first round of experiments. This procedure is presented for the loop fingerprint only. The measured RMSD in these experiments serve as a basis for comparison with the third reconstruction example, where the data-points were picked to be meaningful.

\subsection{Results}\label{subsec:results}
\begin{figure*}
\centering
\subfloat[Fingerprint]{\frame{\includegraphics[width = 0.22\textwidth,trim=6.5cm 9.1cm 5.8cm 8.6cm,clip]{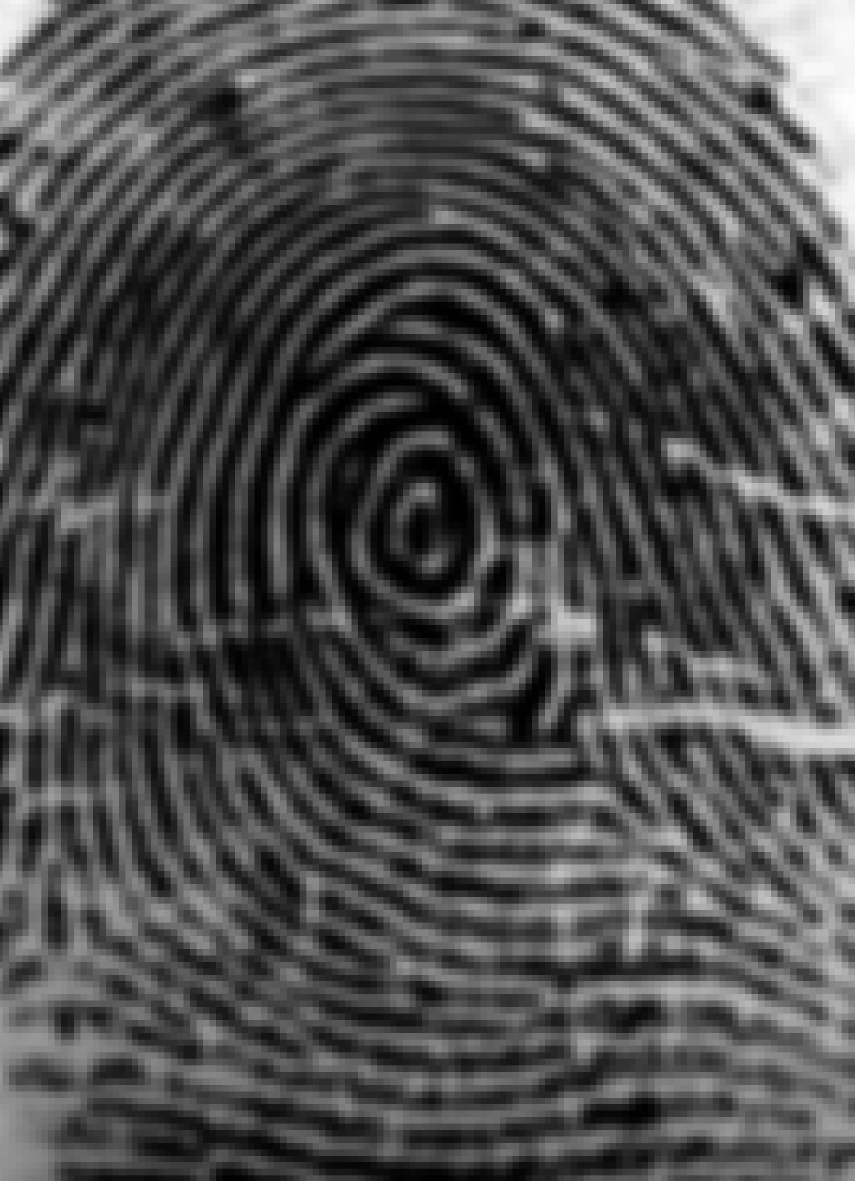}}} 
\quad
\subfloat[Output]{\frame{\includegraphics[width = 0.22\textwidth,trim=6.5cm 9.1cm 5.8cm 8.6cm,clip]{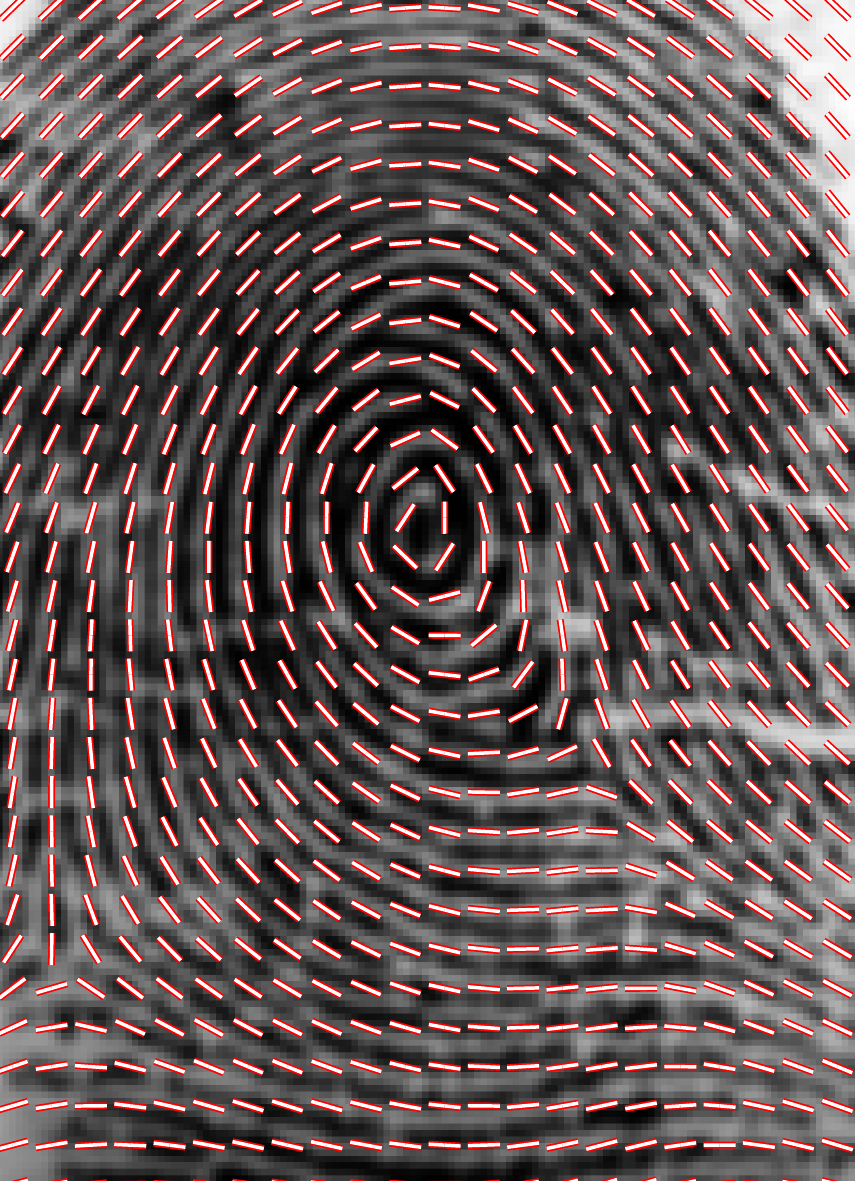}}} 
\quad
\subfloat[Matrix $\overline{\M}$]{\frame{\includegraphics[width = 0.22\textwidth,trim=6.5cm 9.1cm 5.8cm 8.6cm,clip]{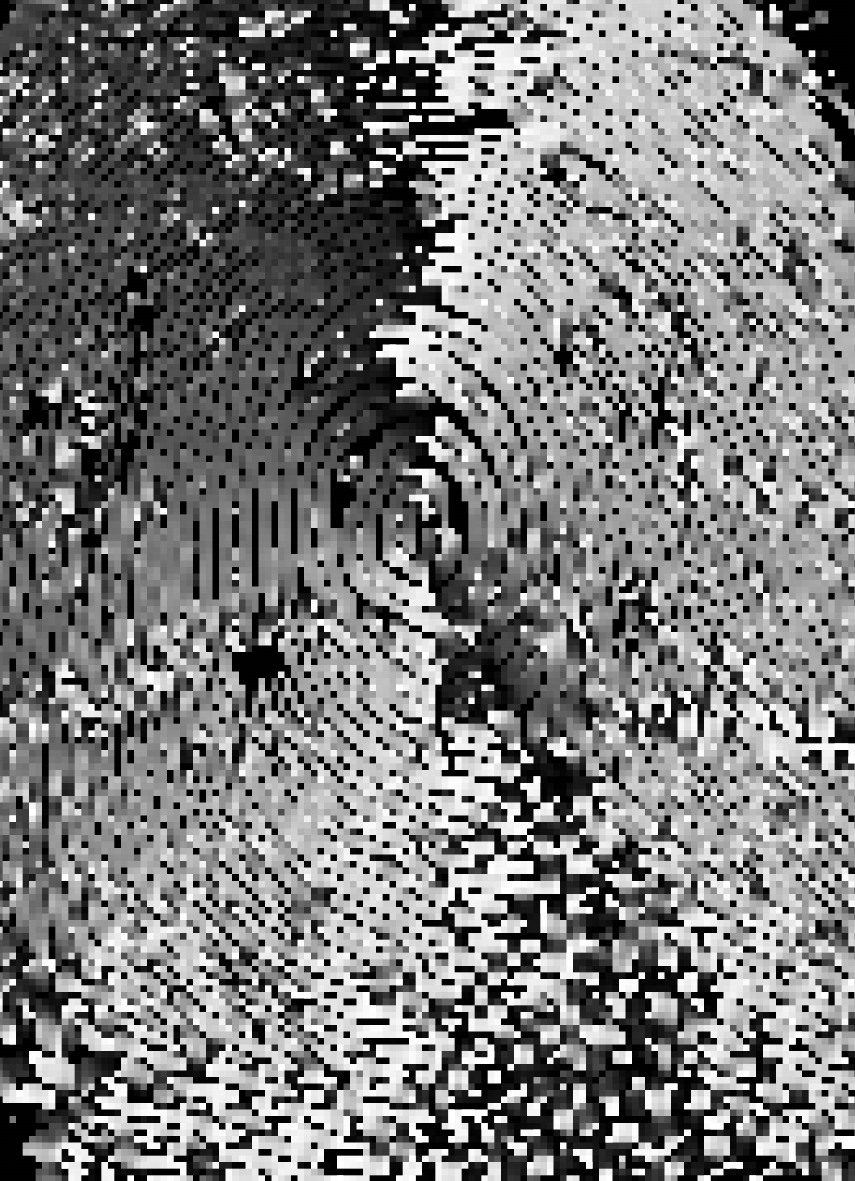}}} 
\quad
\subfloat[Matrix $\overline{\L}$]{\frame{\includegraphics[width = 0.22\textwidth,trim=6.5cm 9.1cm 5.8cm 8.6cm,clip]{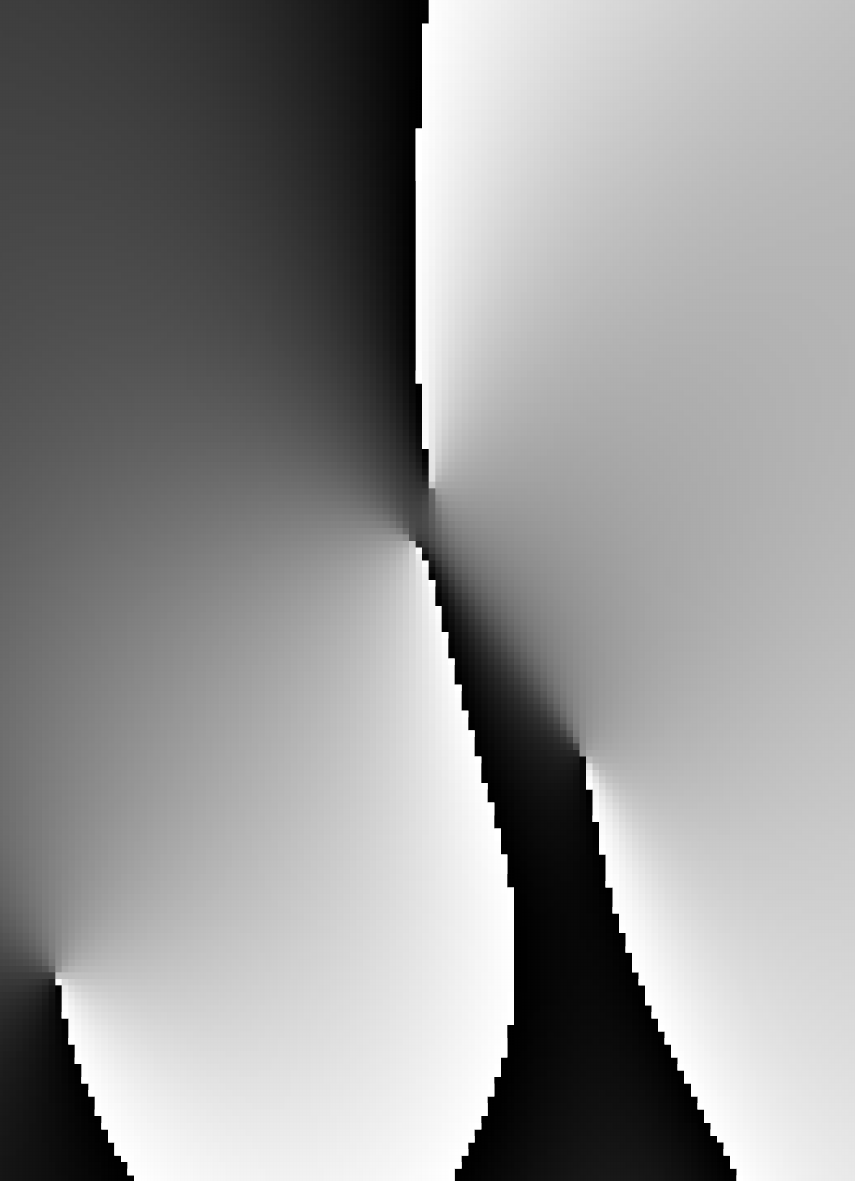}}} 
\quad
\caption{``Whorl'' Experiment ---see details in Section~\ref{subsec:results}}
\label{fig:Whorl-FatSet}
\end{figure*}

\begin{figure*}
\centering
\subfloat[Fingerprint]{\frame{\includegraphics[width = 0.22\textwidth,trim=7.4cm 9.1cm 6.6cm 8.6cm,clip]{Loop_1}}} 
\quad
\subfloat[Output]{\frame{\includegraphics[width = 0.22\textwidth,trim=7.4cm 9.1cm 6.6cm 8.6cm,clip]{Loop_2}}} 
\quad
\subfloat[Matrix $\overline{\M}$]{\frame{\includegraphics[width = 0.22\textwidth,trim=7.4cm 9.1cm 6.6cm 8.6cm,clip]{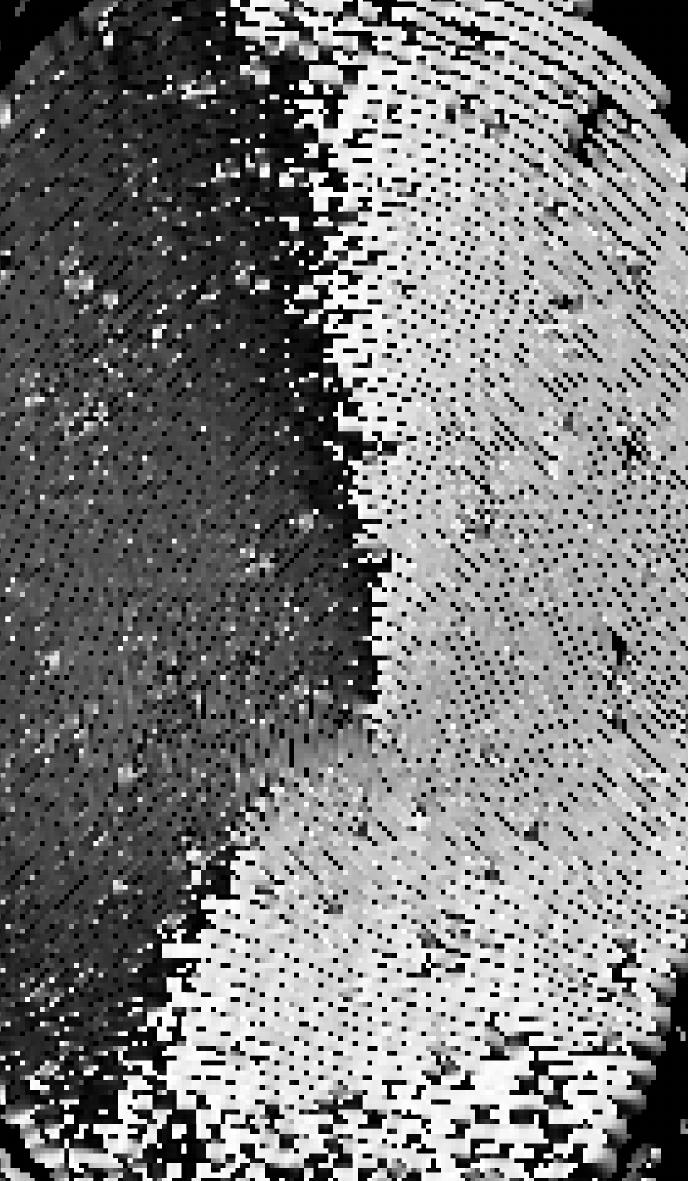}}}
\quad
\subfloat[Matrix $\overline{\L}$]{\frame{\includegraphics[width = 0.22\textwidth,trim=7.4cm 9.1cm 6.6cm 8.6cm,clip]{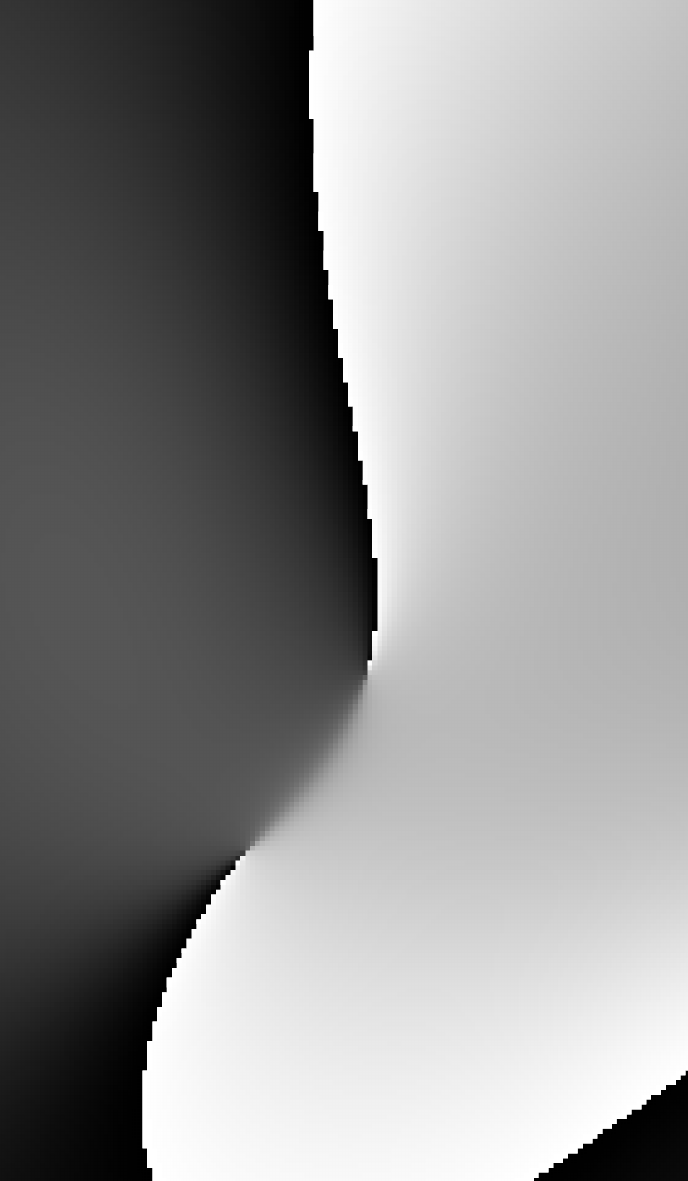}}} 
\quad
\caption{``Loop'' Experiment  ---see details in Section~\ref{subsec:results}}
\label{fig:Loop-FatSet}
\end{figure*}

The result of the interpolation experiments with large datasets are presented in Figures~\ref{fig:Whorl-FatSet}~and~\ref{fig:Loop-FatSet}.
The first image is the input of the algorithm, a grey-scale image of a fingerprint. The second is a representation of  the reconstructed orientation field, as a field of lines, that has been superimposed on the input.

The third and fourth images correspond to a representation of the phases of the two orientation fields, that is, grey-scale images where angles from $0$ to $\pi$ are mapped to a light intensity (near $0$, dark, near $\pi$, bright). Image (c) is then a representation of the  matrix of targeted orientations $\overline{\M}=\left(\M(\bar x_j,\bar y_j)\right)_{j \in \mathcal{J}}$, while image (d) is the matrix $\overline{\L}=\left(\L(\bar x_j,\bar y_j)\right)_{j \in \mathcal{J}}$ of the interpolated field.
Notice that the hard lines separating black and white correspond to the location of transitions from $0$ to $\pi$ in the orientation. Likewise, singularities are points where all grey levels accumulate (similar to the pinwheel singularities observed in the visual cortex \cite{pet2017}).

\begin{figure*}
\centering
\subfloat[Fingerprint]{\frame{\includegraphics[width = 0.22\textwidth,trim=7.4cm 9.1cm 6.6cm 8.6cm,clip]{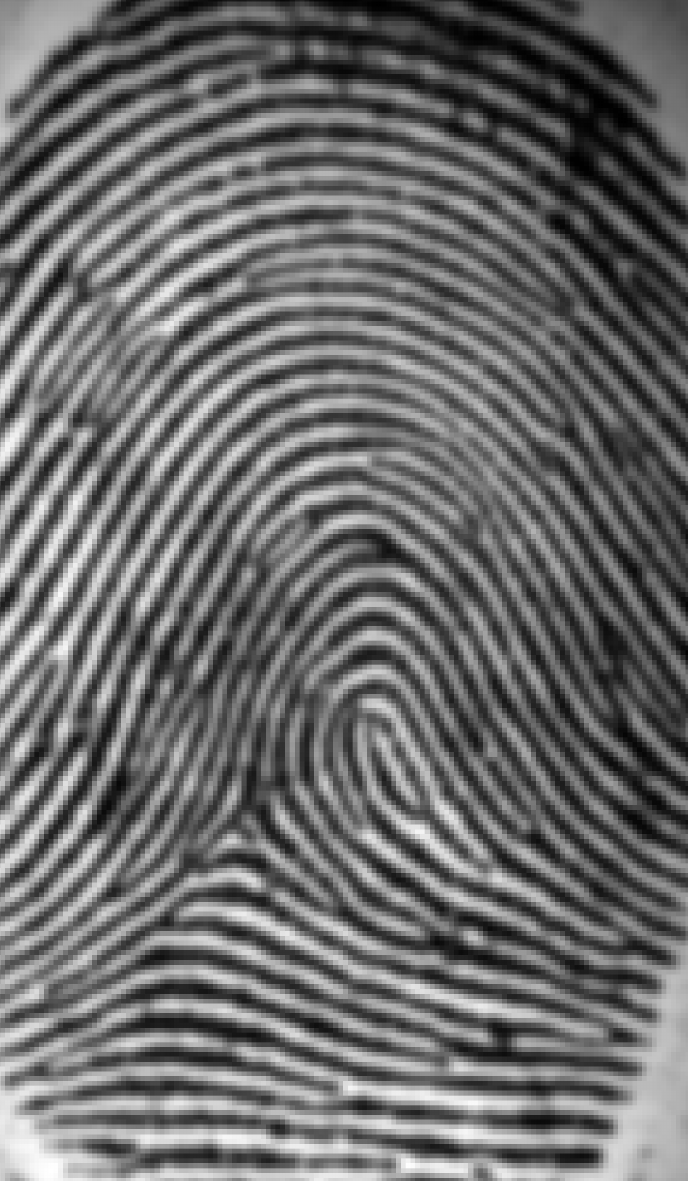}}} 
\quad
\subfloat[Input data]{\frame{\includegraphics[width = 0.22\textwidth,trim=7.4cm 9.1cm 6.6cm 8.6cm,clip]{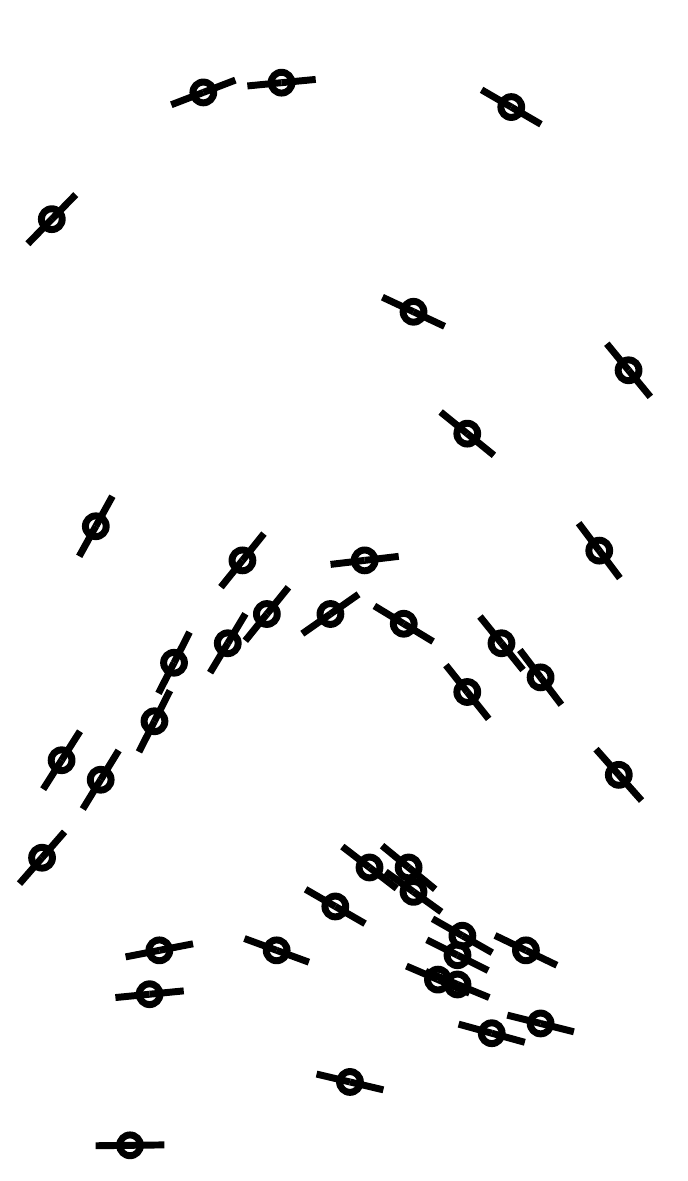}\label{Fig:ScarceInput}}} 
\quad
\subfloat[Output]{\frame{\includegraphics[width = 0.22\textwidth,trim=7.4cm 9.1cm 6.6cm 8.6cm,clip]{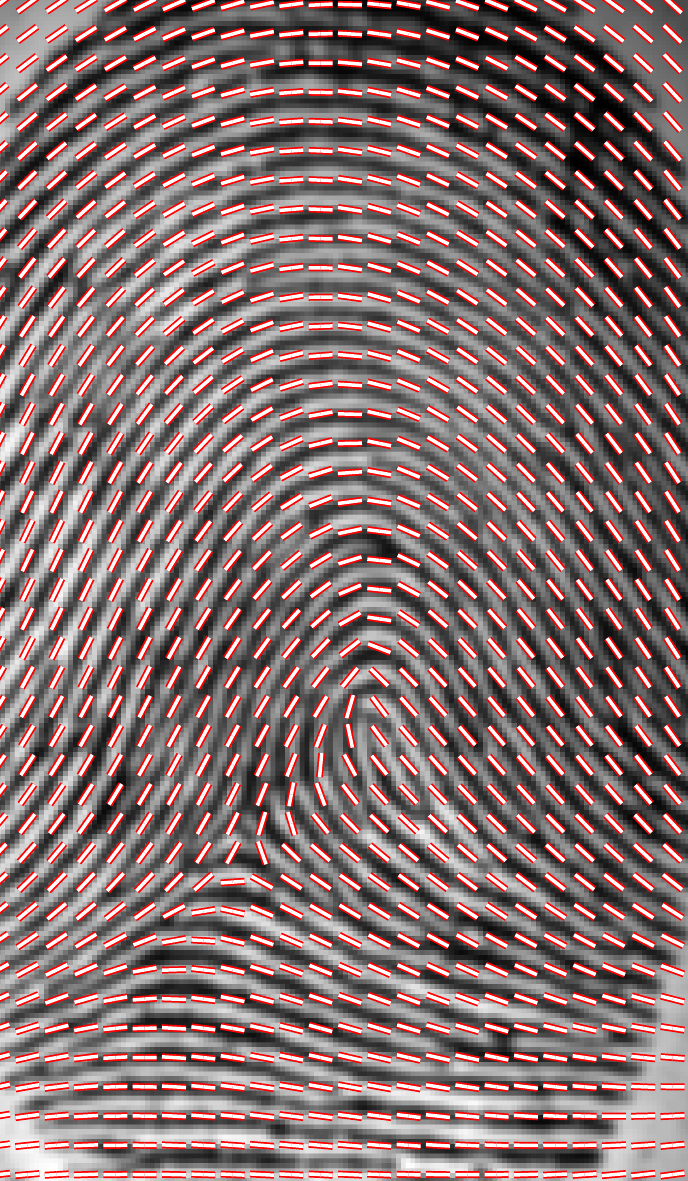}}}
\quad
\subfloat[Matrix $\overline{\L}$]{\frame{\includegraphics[width = 0.22\textwidth,trim=7.4cm 9.1cm 6.6cm 8.6cm,clip]{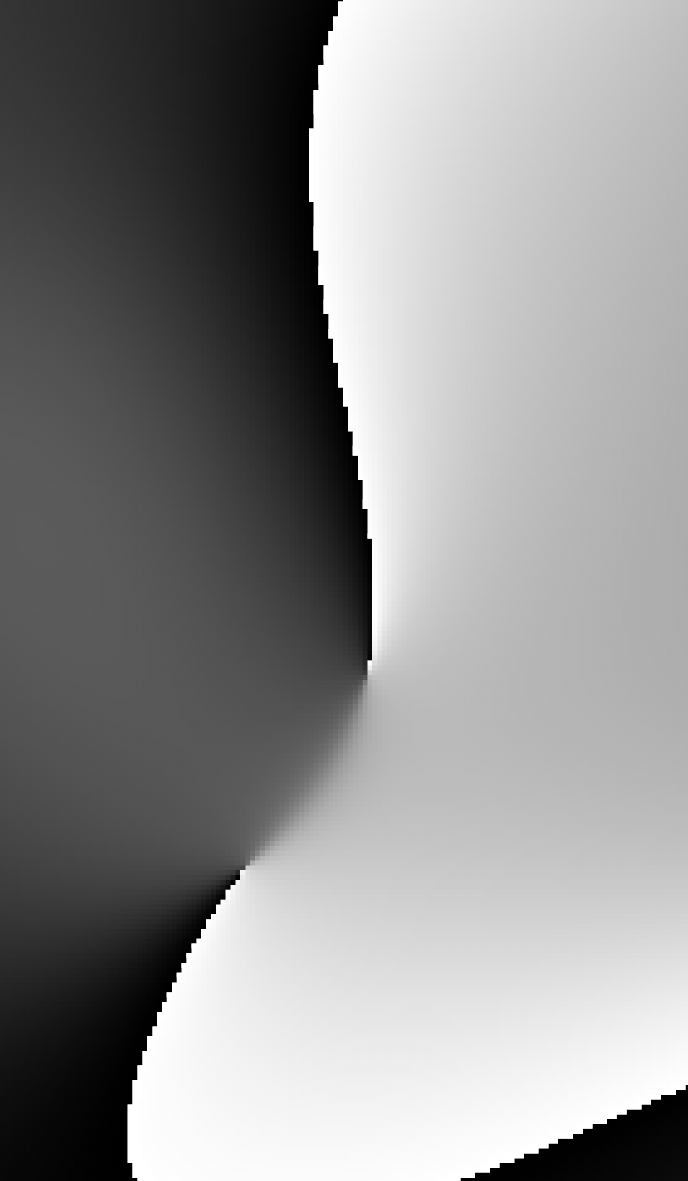}}}
\quad
\caption{``Loop'' Experiment with a dataset made of 40 selected elements  ---see details in Section~\ref{subsec:results}}
\label{fig:Loop-SparseSet}
\end{figure*}

In the case of a scarce dataset, the input is a small collection of triples $(x,y,\theta)\in \R^2\times \P^1$, hence, the results are presented in  Figure~\ref{fig:Loop-SparseSet} in a slightly different way. Indeed, the representation of the matrix $\overline{\M}$ is no longer appropriate: in the third experiment, the data-set $(x_i,y_i)_{i\in \I}$ does not coincide with $(\bar x_j,\bar y_j)_{j \in \mathcal{J}}$. Therefore, we added a representation of the $40$ inputs as line elements on the plane, cf.~Figure~\ref{Fig:ScarceInput}.

\begin{figure}
    \centering
    \begin{minipage}{.3\linewidth}
    \subfloat[Whorl details\label{sfig:whorl_detail1}]{
        \frame{\includegraphics[width = \linewidth,trim=8cm 12.5cm 7cm 11.2cm,clip]{Whorl_1.pdf}}}
    \end{minipage}
    \qquad
    \begin{minipage}{.3\linewidth}
    \subfloat[Interpolation\label{sfig:whorl_detail2}]{
        \frame{\includegraphics[width = \linewidth,trim=8cm 12.5cm 7cm 11.2cm,clip]{Whorl_2.pdf}}}
    \end{minipage}
    
    \begin{minipage}{.3\linewidth}
    \subfloat[Loop details\label{sfig:loop_detail1}]{
        \frame{\includegraphics[width = \linewidth,trim=8.3cm 10.7cm 8.2cm 14.5cm,clip]{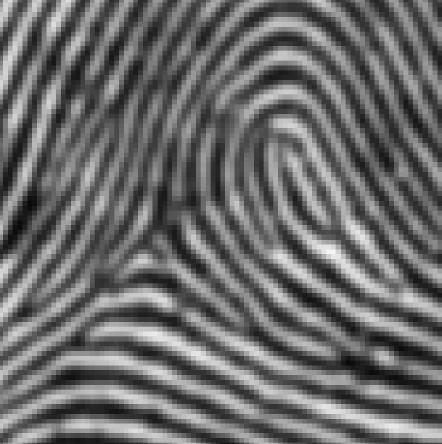}}}
    \end{minipage}
    \qquad
    \begin{minipage}{.3\linewidth}
    \subfloat[Interpolation\label{sfig:loop_detail2}]{
        \frame{\includegraphics[width = \linewidth,trim=8.3cm 10.7cm 8.2cm 14.5cm,clip]{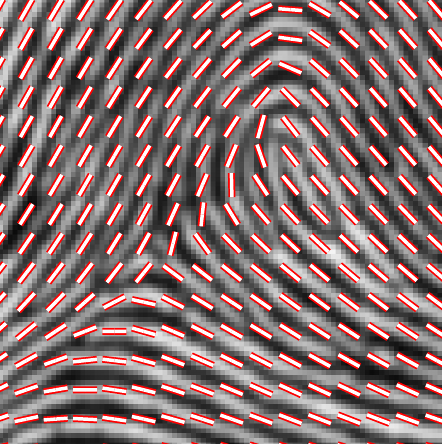}}}
    \end{minipage}

    \begin{minipage}[t]{.3\linewidth}
    \subfloat[Region of inflexion\label{sfig:inflexion_detail1}]{
        \frame{\includegraphics[width = \linewidth,trim=8cm 12.cm 9.5cm 14.2cm,clip]{Loop_1.pdf}}}
    \end{minipage}
    \qquad
    \begin{minipage}[t]{.3\linewidth}
    \subfloat[Imprecise interpolation\label{sfig:inflexion_detail2}]{
        \frame{\includegraphics[width = \linewidth,trim=8cm 12cm 9.5cm 14.2cm,clip]{Loop_2.pdf}}}
    \end{minipage}

    \caption{Details of interpolation figures.}
    \label{fig:details}
    
\end{figure}

Let us briefly comment on these results. As expected, singularities are elegantly fitted without prior knowledge of their positions { (see, for instance, the zoomed in details of Figures~\ref{sfig:whorl_detail1}-\ref{sfig:whorl_detail2} and Figures~\ref{sfig:loop_detail1}-\ref{sfig:loop_detail2})}. Interpolations of degree $(3,3)$ seem to be sufficient to obtain these results and higher order does not improve the final energy significantly. This is in line with other similar interpolation methods (see, for instance, the discussion on the matter in \cite{Wan-Hu2007a}).
It also appears that inflexions in the orientation fields seem to be hard to fit properly 
{ (such as some areas of the Loop experiment, Figure~\ref{fig:Loop-FatSet}, in particular the highlighted region in Figures~\ref{sfig:inflexion_detail1}-\ref{sfig:inflexion_detail2}).}
This is a weakness of the method that seems to be a general difficulty  observed across the different interpolation methods we  encountered in the literature. However, this does not seem to be an obstacle for the later use of finer analysis and reconstruction methods.

In the case of a scarcer dataset, a moderate number of  meaningful points is sufficient to recover some global information on the structure of the orientation field, such as index and curvature of certain regions. However, precise placement of the singularities is dependent on the position of elements chosen as inputs.
Indeed, in the  example case shown in Figure~\ref{fig:Loop-SparseSet},
a majority of data points tightly fit the singular regions which allows a good reconstruction of the underlying orientation field.

\begin{figure*}
\begin{center}

\subfloat[10 data points]{\frame{\includegraphics[width = 0.22\textwidth, trim=7.4cm 9.1cm 6.5cm 8.5cm,clip]{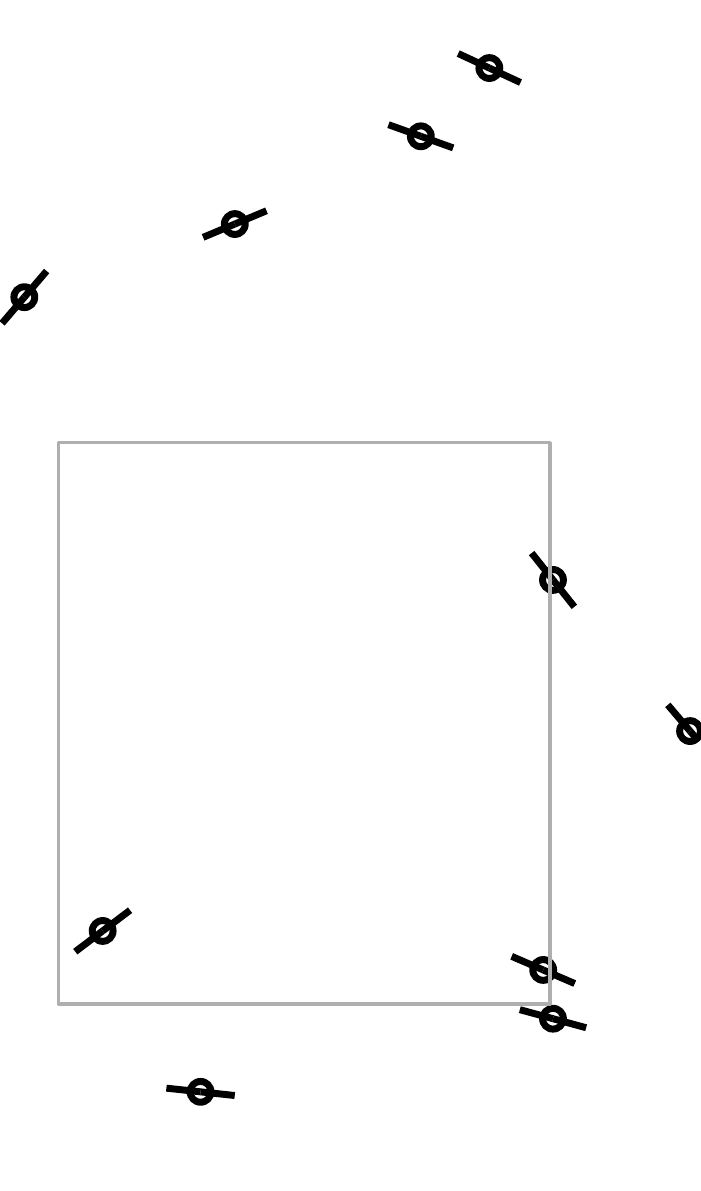}}} 
\quad
\subfloat[20 data points]{\frame{\includegraphics[width = 0.22\textwidth, trim=7.4cm 9.1cm 6.5cm 8.5cm,clip]{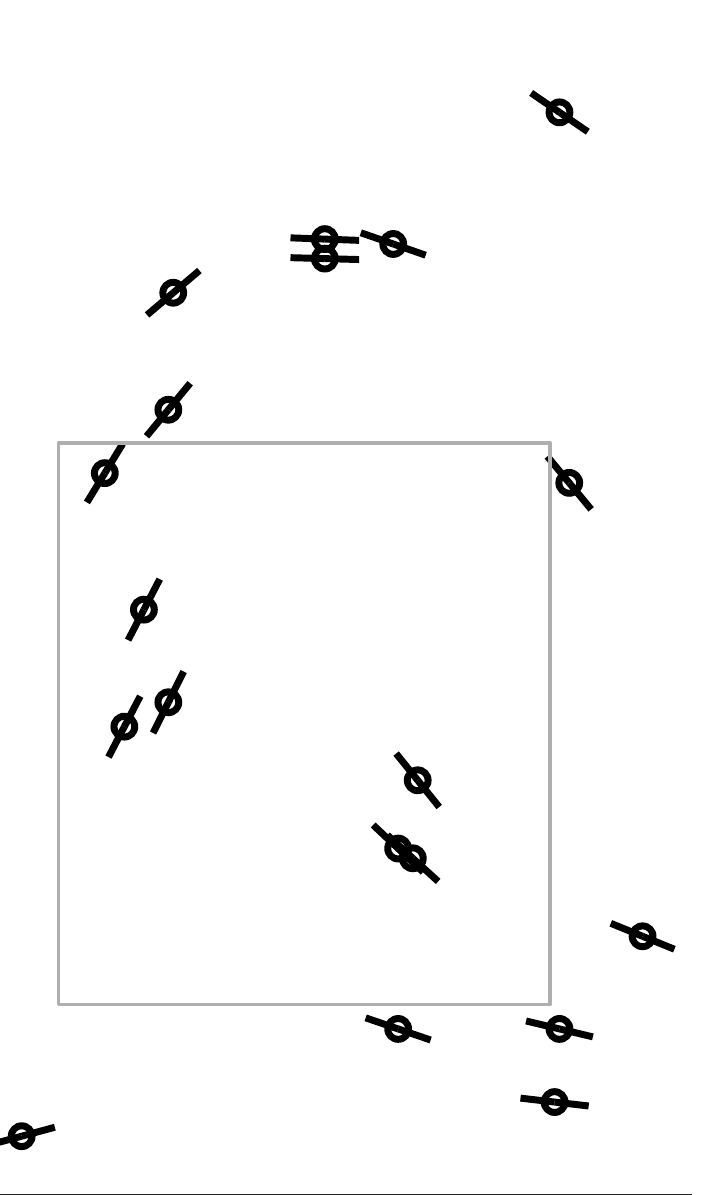}}} 
\quad
\subfloat[40 data points]{\frame{\includegraphics[width = 0.22\textwidth, trim=7.4cm 9.1cm 6.5cm 8.5cm,clip]{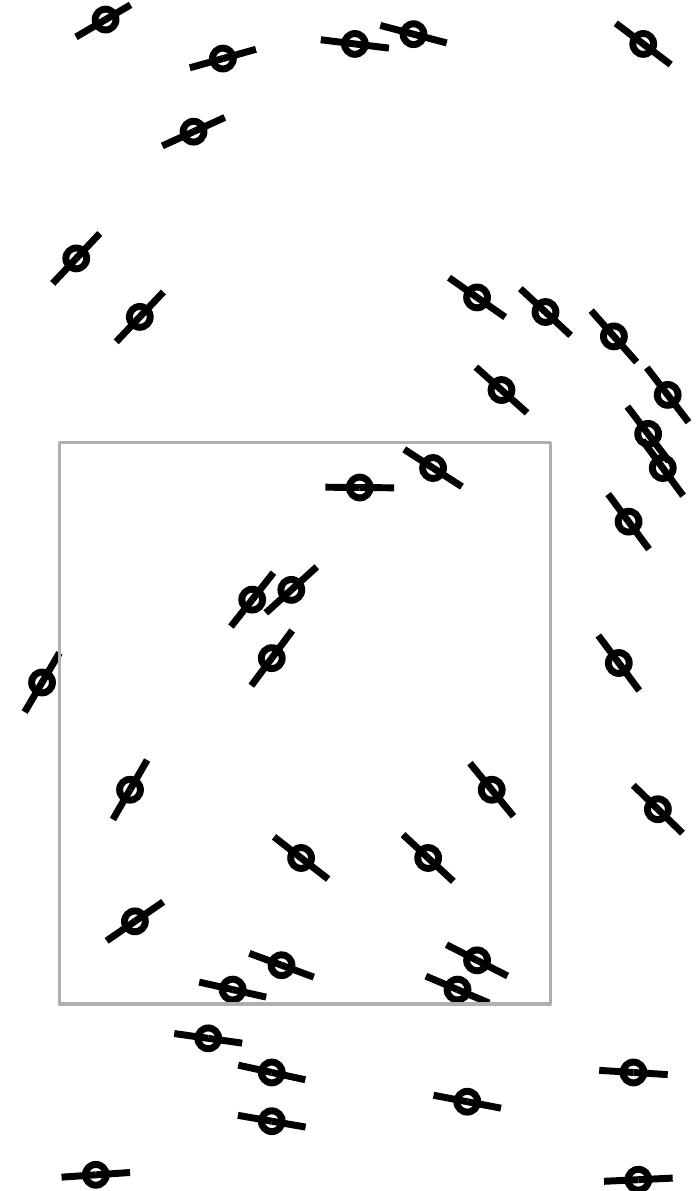}}} 
\quad
\subfloat[80 data points]{\frame{\includegraphics[width = 0.22\textwidth, trim=7.4cm 9.1cm 6.5cm 8.5cm,clip]{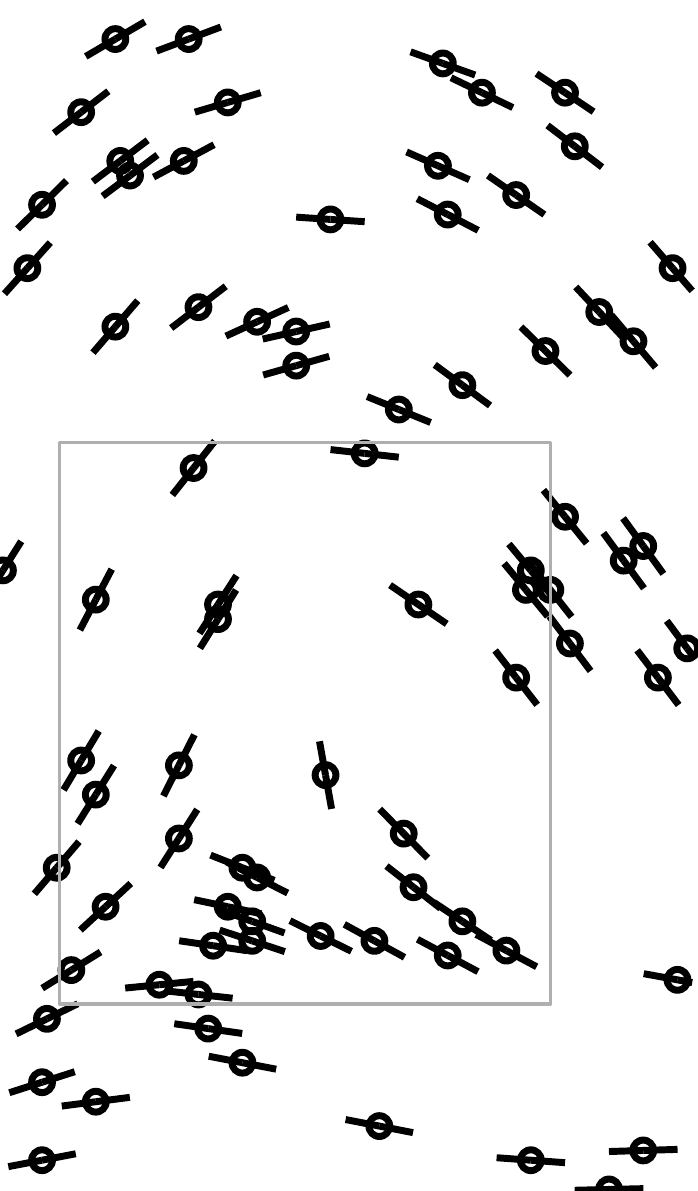}}} 

\subfloat[RMSD = $3.366\times 10^{-1}$\label{fig:SCARCE10_C} ]{\frame{\includegraphics[width = 0.22\linewidth, trim=8.cm 11cm 8.cm 13cm, clip]{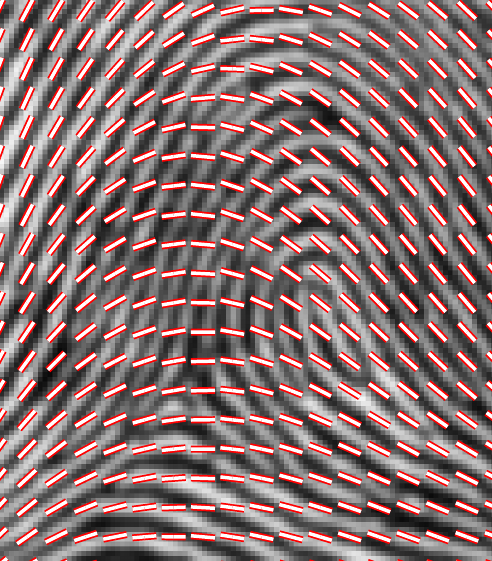}}} 
\quad
\subfloat[RMSD = $2.893\times 10^{-1}$\label{fig:SCARCE20_C}  ]{\frame{\includegraphics[width = 0.22\textwidth, trim=8cm 11cm 8cm 13cm, clip]{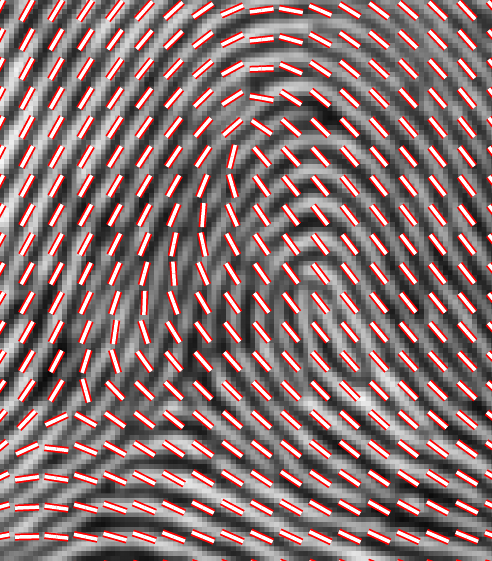}}}
\quad
\subfloat[RMSD = $1.301 \times 10^{-1}$ \label{fig:SCARCE40_C}]{\frame{\includegraphics[width = 0.22\textwidth, trim=8cm 11cm 8cm 13cm, clip]{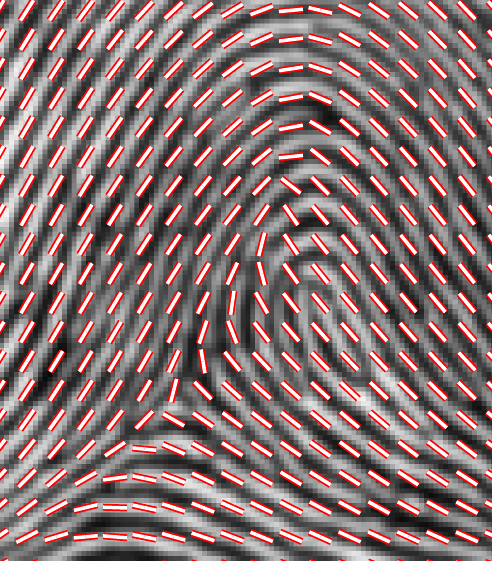}}} 
\quad
\subfloat[RMSD = $0.744\times 10^{-1}$ \label{fig:SCARCE80_C}]{\frame{\includegraphics[width = 0.22\textwidth, trim=8cm 11cm 8cm 13cm, clip]{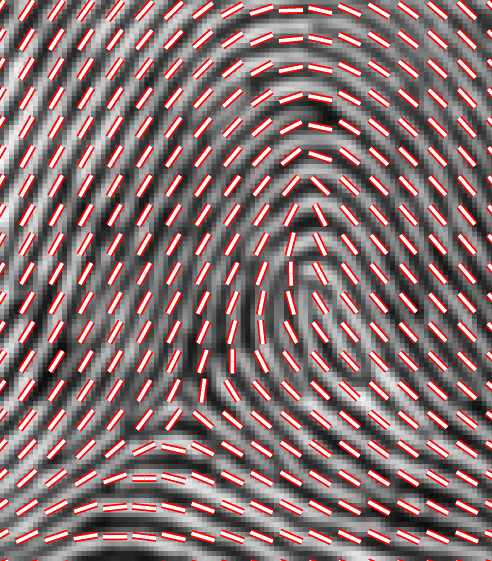}}}

\caption{``Loop” Experiment with typical randomly generated small datasets. The first line represents the input datasets. The corresponding outputs are shown on the second line.
In the outputs, we focus on the delimited regions in order to display the improvements of the interpolation.}
\label{fig:Loop-SparseSet-Random}
\end{center}
\end{figure*}

In contrast, we ran a series of experiments with data-points positions randomly picked in the domain and corresponding orientations taken from the results of the first two experiments. As such, these are considered to be high-fidelity data in the sense that the first reconstruction filtered out the noise. Four series of experiments have been performed with $10$, $20$, $40$ and $80$ data-points. Each batch was made of $10$ datasets. The average RMSD and the standard deviation for each category are reported in Table~\ref{tab:SCARCE_RMSD}. { Representative examples of the batches are shown in Figure~\ref{fig:Loop-SparseSet-Random}.
With only 10 data-points, relevant information is easily missed, and singularities are not reconstructed, see Figure~\ref{fig:SCARCE10_C}.  The singularities appear with more datapoints, but their distribution influences the quality of the placement, see Figure~\ref{fig:SCARCE20_C}-\ref{fig:SCARCE40_C}-\ref{fig:SCARCE80_C}.
} 

As one can expect, the mean RMSD decreases with an increasing number of data-points, and the standard deviation narrows down. However, the RMSD of the third experiment is 0.07278 (experiment corresponding to Figure~\ref{fig:Loop-SparseSet}).
This value, obtained with 40 hand-picked data-points, turns out to be more in line with the values observed with 80 points, highlighting the influence of the singularities on the global pattern.

\begin{table}[tb]
    \centering
    \begin{tabular}{| c|c | c|}
    \hline
           &  RMSD mean  & RMSD standard deviation  \\
           \hline
           \hline
        $\left \vert \I  \right \vert = 10 $ &  $3.809  \times 10^{-1}$ & $1.121 \times 10^{-1}$  \\
 $\left \vert \I  \right \vert = 20 $ &  $2.637 \times 10^{-1} $ & $8.150 \times 10^{-2}$ \\
  $\left \vert \I  \right \vert = 40 $ & $1.519 \times 10^{-1}$ &  $5.978 \times 10^{-2}$ \\
   $\left \vert \I  \right \vert = 80 $ & $ 7.873\times 10^{-2}$  & $2.346 \times 10^{-2}$
   \\ \hline
    \end{tabular}
    \caption{Distribution of final RMSD.}
    \label{tab:SCARCE_RMSD}
\end{table}


%
%

\subsection{Discussion}

As illustrated in the simulations,
the algorithm we designed is satisfactory.
Through a rudimentary gradient based method, it is possible to recover a smooth interpolation of the targeted orientation field while discarding some of the shortcomings of prior methods. For instance the continuity issue is solved with the introduction of the bisector model and a proper energy functional. 
Indeed, the doubling phase step is based on extracting a pair of discontinuous real valued functions which are, in a second step, interpolated with smooth functions.
Conversely, we were able to provide a procedure that performs a smooth interpolation with the true orientation field as direct target.

Furthermore, let us mention some of the applications of the bisector line field interpolation in the context of the study of orientation fields. This generalization of the doubling phase with smooth functions allows the study of the line field with methods geared towards smooth functions. 
Regarding the study of singularities, it is well known that interpolation can be used to locate and describe singularities of the orientation field. In the case of bisector line fields, this observation still holds true and we can illustrate it with an example. Singularities of $\B(\X,\Y)$ correspond to points of $\R^2$ such that either $\X$ or $\Y$ vanish. For instance, if $p\in \R^2$ is such that $\X(p)=0$ and $\Y(p)\neq 0$, the bisector line field $q\mapsto \B(\mathrm{D}_p\X\cdot q,\Y(p))$ acts as a linearization of $\B(\X,\Y)$
 at $p$ as soon as $\mathrm{D}_p\X$ is invertible (which can be assumed as this is the generic case). As a consequence, for instance,  one has  in this case 
$$
\ind_p \B(\X,\Y)=\frac{1}{2}\mathrm{sign}(\det \mathrm{D}_p \X).
$$

The energy we proposed is the  natural choice when considering this problem and the gradient descent method is one possible direction to optimize it. Furthermore, its definition is adaptable and more can be done when information on the target is known.  For instance, one can introduce  weights depending on the quality of the first lift or prior knowledge of the singularity locations.

Finally, solving this optimization problem on the set of polynomial vector fields of a fixed degree is not a requirement of the method.
What is actually necessary to solve the problem with this methodology is the introduction of a family of smooth functions that serves as a basis for the space of regular vector fields on a bounded domain of $\R^2$. Hence this entire method can be adapted to the classical case of trigonometric polynomials.

\section{Conclusion}
In this paper we proposed a solution to the problem of interpolation of orientation fields with smooth functions.
To this end, we introduced a methodology based on the bisector line field model associated with a well suited energy functional. On the one hand, the bisector model has the double virtue of generalizing known techniques from the field of fingerprint analysis while resolving the continuity issues from the classical approach. On the second hand, the energy is coherent with unique aspects of this problem on the space of orientation fields and facilitates to use of gradient descent methods. Finally, the procedure has been applied to perform polynomial interpolation of orientation fields in the framework of fingerprint analysis.

\bibliographystyle{plain}      
\bibliography{OFbib}   

\end{document}